\documentclass[12pt]{amsart}
\usepackage{amsfonts}
\usepackage{amssymb}
\usepackage{amsmath,amscd}
\input xy
\xyoption {all}

\usepackage{geometry}
\usepackage{amsthm}
\usepackage{amssymb}
\usepackage{latexsym}
\usepackage[dvips]{graphics}

\newtheorem{thm}{Theorem}[section] 
\newtheorem{cor}[thm]{Corollary}
\newtheorem{defn}[thm]{Definition}
\newtheorem{example}[thm]{Example}

\newtheorem{lemma}[thm]{Lemma}
\newtheorem{prop}[thm]{Proposition}
\newtheorem{remark}[thm]{Remark}

\newcommand{\DC}{{\mathcal D}}
\newcommand{\LL}{{\mathcal L}}

\newcommand{\M}{{\mathcal M}}

\newcommand{\OO}{{\mathcal O}}
\newcommand{\HC}{{\mathcal H}}
\newcommand{\VV}{{\mathcal V}}
\newcommand{\FC}{{\mathcal F}}

\newcommand{\Z}{\mathbb{Z}}

\newcommand{\Q}{\mathbb{Q}}
\newcommand{\A}{\mathbb{A}}
\newcommand{\Lef}{\mathbb{L}}
\newcommand{\R}{\mathbb{R}}
\newcommand{\C}{\mathbb{C}}

\newcommand{\VB}{\mathbb{V}}

\begin{document}

\title{Characteristic classes of mixed Hodge modules}

\author[J\"org Sch\"urmann ]{J\"org Sch\"urmann}
\address{Dr. J.  Sch\"urmann : Mathematische Institut,
          Universit\"at M\"unster,
          Einsteinstr. 62, 48149 M\"unster,
          Germany.}
\email {jschuerm@math.uni-muenster.de}

\begin{abstract}
This paper is an extended version of an expository  talk  given at the workshop
``Topology of stratified spaces'' at MSRI Berkeley in September 2008.
It gives an introduction and overview about recent developments on
the interaction of the theories of characteristic classes and mixed Hodge
theory for singular spaces in the complex algebraic context.
      It uses M. Saito's deep theory of mixed Hodge modules as a``black
box", thinking about them as ``constructible or perverse sheaves of
Hodge structures", having the same functorial calculus of
Grothendieck functors. For the ``constant Hodge sheaf", one gets the
``motivic characteristic classes" of Brasselet-Sch\"urmann-Yokura,
whereas the classes of the ``intersection homology Hodge sheaf" were
studied by Cappell-Maxim-Shaneson. The classes associated to 
``good"  variation of mixed Hodge structures where studied in connection
with understanding the monodromy action by Cappell-Libgober-Maxim-Shaneson
and the author. 
      There are two versions of these characteristic classes. The
K-theoretical classes capture information about the graded pieces of
the filtered de Rham complex of the filtered D-module underlying a
mixed Hodge module. Application of a suitable Todd class
transformation then gives classes in homology. These classes are
functorial for proper pushdown and exterior products, together with
some other properties one would expect for a ``good" theory of
characteristic classes for singular spaces.
      For ``good"  variation of mixed Hodge structures they have an
explicit classical description in terms of ``logarithmic de Rham
complexes". On a point space they correspond to a specialization of
the Hodge polynomial of a mixed Hodge structure, which one gets by
forgetting the weight filtration.
    Finally also some relations to other subjects of the conference,
like index theorems, signature, L-classes, elliptic genera and motivic characteristic classes for singular spaces,
will be indicated.
\end{abstract}

\maketitle

\tableofcontents

\section{Introduction}
This paper gives an introduction and overview about recent developments on
the interaction of the theories of characteristic classes and mixed Hodge
theory for singular spaces in the complex algebraic context. The reader is not assumed to
have any background on one of these subjects, and the paper can also be used
as a bridge for communication between researchers on one of these subjects.\\

General references for the theory of characteristic classes of singular spaces is the survey
\cite{SY} as well as the paper \cite{Y} in these proceedings. As references for  mixed Hodge
theory one can use \cite{PS, Voi}, as well as the nice paper \cite{P} for explaining
the motivic viewpoint to mixed Hodge theory. Finally as an introduction to  M. Saito's deep theory of mixed Hodge modules
one can use \cite{PS}[chap. 14], \cite{Sa2} as well as the introduction \cite{Sab}.
The  theory of mixed Hodge modules is used here more or less as a``black
box", thinking about them as ``constructible or perverse sheaves of
Hodge structures", having the same functorial calculus of
Grothendieck functors.  The underlying theory of  constructible and perverse sheaves can be found in
\cite{BBD, KS, Sc}.\\

For the ``constant Hodge sheaf" $\Q^H_Z$ one gets the
``motivic characteristic classes" of Brasselet-Sch\"urmann-Yokura \cite{BSY} as explained in \cite{Y} in these proceedings.
The classes of the ``intersection homology Hodge sheaf" $IC_Z^H$ were
studied by Cappell-Maxim-Shaneson in \cite{CMS0, CMS}.
Also, the classes associated to ``good"  variation of mixed Hodge structures  where studied via Atiyah-Meyer type formulae
by Cappell-Libgober-Maxim-Shaneson in \cite{CLMS, CLMS2}.
For a summary compare also with \cite{MSc}.\\

There are two versions of these characteristic classes, the {\em motivic Chern class transformation} $MHC_y$ and the {\em motivic
Hirzebruch class transformation} $MHT_{y*}$.
The $K$-theoretical classes $MHC_y$ capture information about the graded pieces of
the filtered de Rham complex of the filtered D-module underlying a
mixed Hodge module. Application of a suitable twisting $td_{(1+y)}$ of the Todd class
transformation $td_*$ of Baum-Fulton-MacPherson \cite{BFM, Fu} then gives the classes 
$MHT_{y*}= td_{(1+y)}\circ MHC_y$ in homology. It is the {\em motivic Hirzebruch class transformation} $MHT_{y*}$, which unifies
\begin{enumerate}
 \item[(-1)] the (rationalized) {\em Chern class transformation} $c_*$ of MacPherson \cite{M}, 
\item[(0)] the {\em Todd class transformation} $td_*$ of
Baum-Fulton-MacPherson \cite{BFM}, and 
\item[(1)] the {\em $L$-class  transformation} $L_*$ 
of Cappell-Shaneson \cite{CS} 
\end{enumerate}
for $y=-1,0$ and $1$ respectively (compare with \cite{BSY, SY} and also with \cite{Y} in these proceedins).
But in this paper we focus on the  $K$-theoretical classes $MHC_y$, because these imply then also the corresponding
results for $MHT_{y*}$ just by application of the (twisted) Todd class transformation.
So the {\em motivic Chern class transformation} $MHC_y$ studied here is really the basic one!\\

Here we explain  the functorial calculus of these classes, stating first in a very precise form the
key results used from Saito's theory of mixed Hodge modules, and explaing then how to get from this the basic results about the
motivic Chern class transformation $MHC_y$. Moreover these results are illustrated by many interesting examples.
For the convenience of the reader, the most general results are only stated in the end of the paper.
In fact, while most of the paper is a detailed survey of the K-theoretical version of the theory as developed in
\cite{BSY, CLMS, CLMS2, MSc}, it is this last section which contains new results on the important functorial properties of these
characteristic classes.
The first two section do not use mixed Hodge modules and are formulated in the (now) classical language of
(variation of) mixed Hodge structures.
Here is the plan of the paper:
\begin{description}
 \item[Section 2] gives an introduction to pure and mixed Hodge structures and the corresponding Hodge genera like
$E$-polynomial and $\chi_y$-genus. These are suitable generating functions of Hodge numbers with $\chi_y$ using only the
Hodge filtration $F$, whereas the $E$-polynomial also uses the weight filtration. We also carefully explain,
why only the $\chi_y$-genus can be further generalized to characteristic classes, i.e. why one has to forget the
weight filtration for applications to characteristic classes.
\item[Section 3] motivates and explains the notion of a variation (or family) of pure and  mixed Hodge structures over a smooth
(or maybe singular) base. Basic examples come from the cohomology of the fibers of a family of complex algebraic varieties.
We also introduce the notion of a ``good'' variation of mixed Hodge structures on a complex algebraic manifold $M$, to shorten the notion for
 a graded polarizable
variation of mixed Hodge structures on $M$, which is {\em admissible} in the sense of Steenbrink-Zucker \cite{SZ} and Kashiwara \cite{Ka},
with {\em quasi-unipotent monodromy} at infinity, i.e. with respect to a compactification $\bar{M}$ of $M$ by a compact complex algebraic manifold $\bar{M}$, with complement $D:=\bar{M}\backslash M$ a normal crossing divisor with smooth irreducible components.
Later on these will give the basic example of so called ``smooth'' mixed Hodge modules.
And for these  good variations we introduce a simple {\em cohomological} characterstic class transformtion $MHC^y$,
which behaves nicely with respect to smooth pullback, duality and (exterior) products. 
As a first approximation to more general mixed Hodge modules and their characteristic classes, we also study in detail
functorial properties of the canonical Deligne extension across a  normal crossing divisor $D$ at infinity (as above),
leading to {\em cohomological} characteristic classes $MHC^y(j_*(\cdot))$ defined in terms of ``logarithmic de Rham
complexes". These classes of good variations have been studied in detail in \cite{CLMS, CLMS2, MSc}, and most results described
here are new functorial reformulations of the results from these sources.
\item[Section 4] starts with an introduction to Saito's functorial theory of algebraic mixed Hodge modules,
explaining its power in many examples, e.g. how to get a pure Hodge structure on the global Intersection
cohomology $IH^*(Z)$ of a compact complex algebraic variety $Z$. From this we deduce the basic calculus of
Grothendieck groups $K_0(MHM(\cdot ))$ of mixed Hodge modules needed for our motivic Chern class transformation $MHC_y$.
We also explain the relation to the motivic view point coming from relative Grothendieck groups of complex
algebraic varieties.
\item[ Section 5.1] is devoted to the definition of our motivic characteristic {\em homology class} transformations $MHC_y$ and $MHT_{y*}$ for
 mixed Hodge modules. By Saito's theory they commute with push down for proper morphisms, and on a compact space
one gets back the corresponding $\chi_y$-genus by pushing down to a point, i.e. by taking the degree of these
characteristic homology classes.
 \item[Sections 5.2-5.3] finally explain other important functoriality properties, like
\begin{enumerate}
 \item Multiplicativity for exterior products.
\item The behaviour under smooth pullback given by a Verdier Riemann-Roch formula.
\item A ``going up and down'' formula for proper smooth morphisms.
\item Multiplicativity between $MHC^y$ and $MHC_y$ for a suitable (co)homological pairing in the context of a morphism with smooth target.
As special cases one gets from this interesting Atiyah and Atiyah-Meyer type formulae
(as studied in \cite{CLMS, CLMS2, MSc}).
\item The relation between $MHC_y$ and duality, i.e. the Grothendieck duality transformation for coherent sheaves ond the Verdier duality for
mixed Hodge modules.
\item The identification of $MHT_{-1*}$ with the (rationalized) Chern class transformation $c_*\otimes\Q$ of MacPherson for the underlying
constructible sheaf complex or function.
\end{enumerate}
\end{description}

Note that such a functorial calculus is expected for any good theory of functorial characteristic classes of singular spaces
(compare \cite{BSY, SY}):
\begin{enumerate}
 \item[c:] For MacPherson's Chern class transformation $c_*$ compare with \cite{BSY, Ke, M,  SY}.
\item[td:] For Baum-Fulton-MacPherson's Todd class transformation $td_*$  compare with \cite{BFM, BFM2, BSY, Fu, FM,  SY}.
\item[L:] For Cappel-Shaneson's $L$-class transformation $L_*$ compare with \cite{BCS, Ba, Ba2, BSY, CS,  SY, Si, Woo}.
\end{enumerate}

Note that the counterpart of mixed Hodge modules in these theories are constructible functions and sheaves (for $c_*$),
coherent sheaves (for $td_*$) and selfdual perverse or constructible sheaf complexes (for $L_*$).
The cohomological counterpart of the smooth mixed Hodge modules (i.e. good variation of mixed Hodge structures) are
locally constant functions and sheaves (for $c^*$), locally free coherent sheaves or vector bundles (for the Chern character $ch^*$) and
selfdual local systems (for $L^*$ and the $KO$-classes of Meyer \cite{Mey}).\\

In this paper we concentrate mainly on pointing out the relation and analogy to the $L$-class story related to important signature invariants,
because these are the subject of many other talks from the conference given in more topological terms.
Finally also some relations to other themes of the conference,
like index theorems, $L^2$-cohomology,  elliptic genera and motivic characteristic classes for singular spaces,
will be indicated.

\section{Hodge structures and genera}
\subsection{Pure Hodge structures}
Let $M$ be a compact {\em K\"{a}hler manifold} (e.g. a complex projective manifold)
of complex dimension $m$.
By classical Hodge theory one gets  the decomposition (for $0\leq n \leq 2m$)
\begin{equation}\label{h-dec}
H^n(M,\C)=\oplus_{p+q=n}\; H^{p,q}(M)
\end{equation}
of the complex cohomology of $M$ into the spaces $H^{p,q}(M)$ of harmonic forms of type
$(p,q)$. This decomposition doesn't depend on the choice of a K\"{a}hler form (or metric) on $M$, and for a complex algebraic manifold $M$ it is of algebraic nature.
Here it is more natural to work with the {\em Hodge filtration}
\begin{equation}\label{h-fil}
F^i(M):=\oplus_{p\geq i}\;H^{p,q}(M)
\end{equation}
so that $ H^{p,q}(M)=F^p(M)\cap \overline{F^q(M)}$, with $\overline{F^q(M)}$ the complex conjugate of $F^q(M)$ with respect to the real structure $H^n(M,\C)=H^n(M,\R)\otimes \C$.
If 
$$ \begin{CD} \Omega_M^{\bullet}=[\OO_M @> d >> \cdots @> d >> \Omega^m_M]
\end{CD}$$ 
denotes the usual holomorphic de Rham complex (with $\OO_M$ in degree zero), then one gets
$$H^*(M,\C)=H^*(M,\Omega_M^{\bullet})$$
by the holomorphic Poincar\'{e}-lemma, and the Hodge filtration is induced from the 
``stupid" decreasing filtration
\begin{equation}\label{dr-fil} \begin{CD}
F^p\Omega_M^{\bullet}=[0@>>> \cdots 0 @>>>  \Omega^p_M @> d >> \cdots @> d >> \Omega^m_M] \:.
\end{CD}
\end{equation}
More precisely, the corresponding {\em Hodge to de Rham spectral-sequence}
degenerates at $E_1$, with 
\begin{equation} \label{Hpq}
E_1^{p,q}=H^{q}(M,\Omega^p_M)\simeq H^{p,q}(M)\:.
\end{equation}
More generally, the same results are true for a compact complex manifold $M$,
which is only {\em bimeromorphic to a  K\"{a}hler manifold} (compare e.g. \cite{PS}[cor.2.30]). This is especially 
true for a compact complex algebraic manifold $M$. Moreover in this case one can calculate
by Serre's GAGA-theorem 
$H^*(M,\Omega_M^{\bullet})$ also with the algebraic (filtered) de Rham complex in the Zariski topology.\\

Abstracting these properties, one can say the $H^n(M,\Q)$ gets an induced {\em pure Hodge
structure of weight $n$} in the following sense:
\begin{defn} Let $V$ be a finite dimesional rational vector space.
A (rational) Hodge
structure of weight $n$ on $V$ is a decomposition
$$V_{\C}:=V\otimes_{\Q}\C = \oplus_{p+q=n}\; V^{p,q}, \quad \text{with $V^{q,p}=\overline{V^{p,q}}$}
\quad \text{(Hodge decomposition).}$$
In terms of the (decreasing) {\em Hodge filtration} $F^iV_{\C}:=\oplus_{p\geq i}\;V^{p,q}$, this is equivalent to the condition
$$F^pV\cap \overline{F^qV}= \{0\} \quad \text{whenever $p+q=n+1$} \quad \text{($n$-opposed filtration).}$$
Then $V^{p,q}=F^p\cap \overline{F^q}$, with $h^{p,q}(V):=dim(V^{p,q})$ the corresponding
{\em Hodge number}.
\end{defn}

If $V,V'$ are rational vector spaces with Hodge structures of weight $n$ and $m$, then
$V\otimes V'$ gets an induced Hodge structure of weight $n+m$, with Hodge filtration
\begin{equation}\label{F-tensor}
F^k(V\otimes V')_{\C}:= \oplus_{i+j=k}\; F^iV_{\C}\otimes F^jV'_{\C} \:.
\end{equation}
Similarly the dual vector space $V^{\vee}$ gets an induced Hodge structure of weight $-n$,
with
\begin{equation}\label{F-dual}
F^k( V^{\vee}_{\C}):=(F^{-k}V_{\C})^{\vee} \:.
\end{equation}
A basic example is the {\em Tate Hodge structure} of weight $-2n\in \Z$ given by
the one dimensional rational vector space
$$\Q(n):=(2\pi i)^n\cdot \Q \subset \C,\quad \text{with $\Q(n)_{\C}=(\Q(n)_{\C})^{-n,-n}$.}$$
Then integration defines an isomorphism $H^2(P^1(\C),\Q)\simeq \Q(-1)$, with
$\Q(-n)=\Q(-1)^{\otimes n}$, $\Q(1)=\Q(-1)^{\vee}$ and $\Q(n)=\Q(1)^{\otimes n}$ for $n>0$.

\begin{defn}\label{pol}
A {\em polarization} of a rational Hodge structure $V$ of weight $n$ is a rational
$(-1)^n$-symmetric bilinear form $S$ on $V$ such that
$$S(F^p,F^{n-p+1})=0 \quad \text{for all $p$ and $i^{p-q}S(u,\bar{u})> 0$ for all
$0\neq u\in V^{p,q}$.}$$
So for $n$ {\em even} one gets in particular
\begin{equation} \label{positive}
(-1)^{p-n/2}S(u,\bar{u})> 0 \quad \text{for all $q$ and
$0\neq u\in V^{p,q}$.}
\end{equation}
$V$ is called {\em polarizable}, if such a polarization exists.
\end{defn}

For example the cohomology $H^n(M,\Q)$ of a projective manifold is polarizable by
the choice of a suitable K\"{a}hler form! Also note that a
polarization of a rational Hodge structure $V$ of weight $n$ induces an isomorphism of
Hodge structures (of weight $n$):
$$V \simeq V^{\vee}(-n):=V^{\vee}\otimes_{\Q}\Q(-n)\:.$$
So if we choose the isomorphism of rational vector spaces $\Q(-n)=(2\pi i)^{-n}\cdot \Q\simeq \Q$,
then a polarisation induces a $(-1)^n$-symmetric duality isomorphism $V\simeq V^{\vee}$.

\subsection{Mixed Hodge structures}
The cohomology (with compact support) $H^n_{(c)}(X,\Q)$ of a singular or non-compact complex algebraic variety can't have a pure Hodge structure in general, but by Deligne 
\cite{De1, De3} it carries a
canonical functorial (graded polarizable) {\em mixed Hodge structure} in the following sense:
\begin{defn}
A finite dimensional rational vector space $V$ has a mixed Hodge structure,
if there is a (finite) increasing {\em weight filtration} $W=W_{\bullet}$ on $V$ (by rational subvector spaces), and a (finite) decreasing Hodge filtration $F=F^{\bullet}$ on $V_{\C}$,
such that $F$ induces a Hodge structure of weight $n$ on $Gr^W_nV:=W_nV/W_{n-1}V$ for all $n$.
Such a mixed Hodge structure is called  {\em (graded) polarizable} if each graded piece 
$Gr^W_nV$ is polarizable.
\end{defn}

A morphism of mixed Hodge structures is just a homomorphism of rational vector spaces
compatible with both filtrations. 
Such a morphism is then {\em strictly} compatible with both filtrations so that
the category $mHs^{(p)}$ of (graded polarizable) mixed Hodge structures is an abelian category,
with  $Gr^W_*, Gr_F^*$ and $Gr_F^*Gr^W_*$ preserving short exact sequences.
$mHs^{(p)}$ is also endowed with a tensor product $\otimes$ and a duality $(\cdot)^{\vee}$, where the corresponding Hodge and weight filtrations are defined as in
(\ref{F-tensor}) and (\ref{F-dual}). So for a complex algebraic variety $X$ one can consider
its cohomology class
$$[H^*_{(c)}(X)]:=\sum_{i}\; (-1)^i�\cdot [H^i_{(c)}(X,\Q)] \in K_0(mHs^{(p)})$$
in the Grothendieck group $K_0(mHs^{(p)})$ of (graded polarizable) mixed Hodge structures.
The functoriality of Deligne's mixed Hodge structure means in particular, that for a
closed complex algebraic subvariety $Y\subset X$, with open complement $U=X\backslash Y$,
the corresponding long exact cohomology sequence
\begin{equation}\label{long-H}
\cdots H^i_{c}(U,\Q)\to H^i_{c}(X,\Q) \to H^i_{c}(Y,\Q)\to \cdots
\end{equation}
is an exact sequence of mixed Hodge structures. Similarly the K\"{u}nneth isomorphism
\begin{equation}\label{K-H}
H^*_{c}(X,\Q)\otimes H^*_{c}(Z,\Q) \simeq H^*_{c}(X\times Z,\Q) 
\end{equation}
for complex algebraic varieties $X,Z$
is an isomorphism of  mixed Hodge structures. Let us denote by $K_0(var/pt)$ the Grothendieck
group of complex algebraic varieties, i.e. the free abelian group of isomorphism classes $[X]$ of such varieties divided out by the {\em additivity relation} 
$$[X]=  [Y]+ [X\backslash Y]$$
for $Y\subset X$ a closed complex subvariety. This is then a commutative ring with addition resp. multiplication induced by the disjoint union resp. the product of varieties.
So by (\ref{long-H}) and (\ref{K-H}) we get an induced ring homomorphism
\begin{equation}\label{H-ring}
 \chi_{Hdg}: K_0(var/pt)\to K_0(mHs^{(p)});\: [X]\mapsto [H^*_{c}(X)] \:.
\end{equation}

\subsection{Hodge genera}
The {\em E-polynomial}
\begin{equation}\label{E}
E(V):=\sum_{p,q}\; h^{p,q}(V)\cdot u^pv^q \in \Z[u^{\pm1},v^{\pm 1}]
\end{equation}
of a rational mixed Hodge structure $V$ with {\em Hodge numbers}
$$h^{p,q}(V):=dim_{\C} Gr_F^pGr^W_{p+q}(V_{\C}) \:,$$
induces a {\em ring} homomorphism
$$E: K_0(mHs^{(p)})\to  \Z[u^{\pm1},v^{\pm 1}]\:,\quad \text{with $E(\Q(-1))=uv$.}$$
Note that $E(V)(u,v)$ is {\em symmetric} in $u$ and $v$, since $h(V)=\sum_n\;h(W_nV)$ and
$V^{q,p}=\overline{V^{p,q}}$ for a pure Hodge structure.
With respect to {\em duality} one has in addition the relation
\begin{equation}\label{E-dual}
E(V^{\vee})(u,v) = E(V)(u^{-1},v^{-1}) \:.
\end{equation}

Later on we will be mainly interested in the following specialized ring homomorphism
$$\chi_y:=E(-y,1): K_0(mHs^{(p)})\to  \Z[y^{\pm1}]\:,\quad \text{with $\chi_y(\Q(-1))=-y$,}$$
defined by
\begin{equation}\label{chi-y-H}
\chi_y(V):= \sum_{p}\;  dim_{\C}(Gr_F^p(V_{\C}))\cdot (-y)^p \:.
\end{equation}
So here one uses only the Hodge and forgets  the weight filtration of a mixed Hodge structure.
With respect to {\em duality} one has then the relation
\begin{equation}\label{chi-y-dual}
\chi_y(V^{\vee}) = \chi_{1/y}(V) \:.
\end{equation}

Note that $\chi_{-1}(V)=dim(V)$ and for a pure polarized Hodge structure $V$ of weight $n$
one has by $\chi_1(V)=(-1)^n\chi_1(V^{\vee})=(-1)^n\chi_1(V)$ and
(\ref{positive}):
$$\chi_1(V)=\begin{cases}
0 &\text{for $n$ odd,}\\
sign(V)&\text{for $n$ even,}
\end{cases}$$
where $sign$ denotes the {\em signature} of the induced symmetric bilinear form $(-1)^{n/2}\cdot S$ on $V$. A similar but deeper result is the famous
{\em Hodge index theorem} (compare e.g. \cite{Voi}[thm.6.3.3])):
$$\chi_1([H^*(M)])=sign(H^m(M,\Q))$$
for $M$ a compact K\"{a}hler manifold of complex even dimension $m=2n$.
Here the right side denotes the signature of the symmetric intersection pairing
$$\begin{CD}
H^m(M,\Q)\times H^m(M,\Q)@> \cup >> H^{2m}(M,\Q)\simeq \Q \:.
\end{CD}$$

The advantage of $\chi_y$ compared to $E$ (and the use of $-y$ in the definition) comes from the following
\begin{description}
\item[Question] Let $E(X):=E([H^*(X)])$ for $X$ a complex algebraic variety. 
For $M$ a compact complex algebraic manifold one gets by (\ref{Hpq}):
$$E(M)=\sum_{p,q\geq 0}\; (-1)^{p+q}\cdot dim_{\C}H^q(M,\Omega_M^p)\cdot u^pv^q \:.$$
Is there a {\em (normalized multiplicative) characteristic class} 
$$cl^*: Iso(\C-VB(M))\to H^*(M)[u^{\pm1},v^{\pm 1}]$$ 
of complex vector bundles such that the E-polynomial is a {\em characteristic number} in the sense that 
\begin{equation}\label{genus}
E(M)=\sharp(M):=deg(cl^*(TM)\cap[M])\in H^*(pt)[u^{\pm1},v^{\pm 1}]
\end{equation}
for any compact complex algebraic manifold $M$ with fundamental class $[M]$?
\end{description}

So the cohomology class $cl^*(V)\in H^*(M)[u^{\pm1},v^{\pm 1}]$ should only depend on the isomorphism class of the complex vector bundle $V$ over $M$ and  commute with pullback.
Multiplicativity says
$$cl^*(V)=cl^*(V')\cup cl^*(V'') \in H^*(M)[u^{\pm1},v^{\pm 1}]$$
for any short exact sequence $0\to V'\to V\to V''\to 0$ of complex vector bundles on $M$.
Finally $cl^*$ is normalized if $cl^*(trivial)=1\in H^*(M)$
for any trivial vector bundle. Then the answer to this question is {\em NO} because
there are unramified coverings $p: M'\to M$ of elliptic curves $M,M'$ of (any) degree $d>0$.
Then $p^*TM\simeq TM'$ and $p_*([M'])=d\cdot [M]$ so that the projection formula
would give for the topological  characteristic numbers the relation 
$$\sharp(M')= d\cdot \sharp(M)\:.$$
But one has
$$E(M)=(1-u)(1-v)=E(M')\neq 0$$
so that the equality $E(M)=\sharp(M)$ is not possible! 
Here wo don't need to ask $cl^*$ to be multiplicative or normalized.
But if we use the invariant $\chi_y(X):=\chi_y([H^*(X)])$, then
$\chi_y(M)=0$ for an elliptic curve, and $\chi_y(M)$ is a characteristic number
in the sense above by the famous {\em generalized Hirzebruch Riemann Roch theorem} (\cite{H}):
\begin{thm}[gHRR]\label{gHRR} There is a unique normalized multiplicative characteristic class 
$$T_y^*: Iso(\C-VB(M))\to H^*(M,\Q)[y]$$ 
such that
$$\chi_y(M)=deg (T_y^*(TM)\cap [M]) = \langle T_y^*(TM), [M]\rangle
\in \Z[y]\subset \Q[y]$$
for any compact complex algebraic manifold $M$. 
Here $\langle \cdot,\cdot\rangle$ is the Kronecker pairing between cohomology and homology.
\end{thm}

The {\em Hirzebruch class} $T^*_y$ and $\chi_y$-genus unify the following (total) characteristic classes and numbers:
$$T^*_y=
\begin{cases}
c^* \:\text{,the Chern class}\\
td^*\:\text{,the Todd class}\\
L^* \:\text{,the  L class}\\
\end{cases} \quad \text{and $\chi_y=$}
\begin{cases}
\chi \:\text{,the  Euler characteristic}\\
\chi_a \:\text{,the arithmetic genus}\\
sign \:\text{,the signature}
\end{cases}
\:\text{for $y=$}
\begin{cases}
-1\\ 0\\1\:. \end{cases}$$

In fact $(gHRR)$ is just a cohomological version of the following $K$-theoretical calculation.
Let $M$ be a compact complex algebraic manifold, so that
\begin{equation} \label{chiy-M} \begin{split}
\chi_y(M)&=\sum_{p,q\geq 0}\; (-1)^{p+q}\cdot  dim_{\C}H^q(M,\Omega_M^p)\cdot (-y)^p \\
&=\sum_{p\geq 0}\; \chi(H^*(M,\Omega_M^p))\cdot y^p \:.
\end{split}\end{equation}
Let us denote by $K^0_{an}(Y)$ (or $G_0^{an}(Y)$)  the Grothendieck group of the exact
(or abelian) category
of holomorphic vector bundles (or coherent $\OO_Y$-module sheaves) on the complex variety $Y$, i.e. the free abelian group
of isomorphism classes $V$ of such vector bundles (or sheaves), divided out by the relation
$$[V]=[V']+[V''] \quad \text{for any short exact sequence $0\to V'\to V\to V''\to 0$.}$$
Then $G_0^{an}(Y)$ (or $K^0_{an}(Y)$) is of (co)homological nature, with
$$f_*: G_0^{an}(X)\to G_0^{an}(Y);\: [\FC]\mapsto \sum_{i\geq 0}\;(-1)^i\cdot [R^if_*\FC]$$
the functorial pushdown for a proper holomorphic map $f: X\to Y$.
In particular for  $X$ compact, the constant map $k: X\to pt$ is proper, with
$$\chi(H^*(X,\FC))=k_*([\FC])\in G^{an}_0(pt)\simeq K^0_{an}(pt)\simeq \Z  \:.$$
Moreover, the tensor product $\otimes_{\OO_Y}$ induces a natural pairing
$$\cap=\otimes: K^0_{an}(Y)\times G_0^{an}(Y) \to G_0^{an}(Y)\:,$$
where we identify a holomorphic vector bundle $V$ with its locally free coherent sheaf of sections $\VV$.
So for $X$ compact we can define a {\em Kronecker pairing}
$$K^0_{an}(X)\times G_0^{an}(X) \to G_0^{an}(pt)\simeq \Z;\:
\langle [\VV],[\FC] \rangle:=k_*([\VV \otimes_{\OO_X} \FC])\:.$$
The {\em total $\lambda$-class} of the dual vector bundle
$$\lambda_y(V^{\vee}):=\sum_{i\geq 0}\; \Lambda^i(V^{\vee})\cdot y^i$$ 
defines a multiplicative characteristic class
$$\lambda_y((\cdot)^{\vee}): K^0_{an}(Y)\to K^0_{an}(Y)[y] \:.$$
And for a compact complex algebraic manifold $M$ one gets the equality
\begin{equation}\label{gHHR-K} \begin{split}
\chi_y(M)&= \sum_{i\geq 0}\; k_*[\Omega^i_M]\cdot y^i\\
&=\langle \lambda_y(T^*M), [\OO_M] \rangle \in G_0^{an}(pt)[y]\simeq \Z[y] \:.
\end{split}\end{equation}

\section{Characteristic classes of variations of mixed Hodge structures}
This section explains the definition of {\em cohomological} characteristic classes associated to good variations of mixed Hodge structures
on complex algebraic and analytic manifolds. These were previously considered in \cite{CLMS, CLMS2, MSc} in connection with Atiyah-Meyer type
formulae of Hodge-theoretic nature. Here we also consider important functorial properties of these classes.
 
\subsection{Variation of Hodge structures}
Let $f: X\to Y$ be a {\em proper smooth} morphism of complex algebraic varieties
or a {\em projective smooth} morphism of complex analytic varieties.
Then the higher direct image sheaf $L=L^n:=R^nf_*\Q_X$ is a {\em locally constant sheaf}
on $Y$ with finite dimensional stalks $$L_y=(R^nf_*\Q_X)_y=H^n(\{f=y\},\Q)$$ 
for $y\in Y$.
Let $\LL:=L\otimes_{\Q_Y} \OO_Y\simeq R^nf_*(\Omega^{\bullet}_{X/Y})$ denote the corresponding holomorphic vector bundle (or locally free sheaf), with $\Omega^{\bullet}_{X/Y}$ the {\em relative holomorphic de Rham complex}.
Then the stupid filtration of $\Omega^{\bullet}_{X/Y}$ determines a decreasing filtration $F$ of $\LL$ by holomorphic subbundles $F^p\LL$, with 
\begin{equation}\label{F-relDR}
Gr_F^p((R^{p+q}f_*\Q_X)\otimes_{\Q_Y} \OO_Y)\simeq R^qf_*(\Omega^{p}_{X/Y})\:,
\end{equation}
inducing for all $y\in Y$ the
Hodge filtration $F$ on the cohomology 
$$H^n(\{f=y\},\Q)\otimes \C\simeq \LL|_y$$
of the compact and smooth algebraic fiber $\{f=y\}$ (compare \cite{PS}[chap.10]). If $Y$ (and therefore also $X$ is smooth),
then $\LL$ gets an induced {\em integrable Gauss-Manin connection} 
$$\nabla: \LL\to \LL\otimes_{\OO_Y}\Omega_Y^1
,\quad \text{with $L\simeq kern(\nabla)$ and $\nabla\circ \nabla=0$,}$$
satisfying the {\em Griffith's transversality} condition
\begin{equation}\label{Griff}
\nabla(F^p\LL)\subset F^{p-1}\LL \otimes_{\OO_Y}\Omega_Y^1 \quad \text{for all $p$.}
\end{equation}
This motivates the following

\begin{defn}
A {\em holomorphic family} $(L,F)$ of Hodge structures of weight $n$ on the reduced complex space $Y$ is a local system $L$ with rational coefficients and finite dimensional stalks on $Y$, and a decreasing filtration $F$ of $\LL=L\otimes_{\Q_Y}\OO_Y$ by holomorphic subbbundles
$F^p\LL$ such that $F$ determines by $L_y\otimes_\Q \C \simeq \LL|_y$
a pure Hodge structure of weight $n$ on each stalk $L_y$ ($y\in Y$).

If $Y$ is a smooth complex manifold, then such a holomorphic family $(L,F)$ is called
a {\em variation} of Hodge structures of weight $n$, if one has in addition
for the induced connection $\nabla: \LL\to \LL\otimes_{\OO_Y}\Omega_Y^1$ the Griffith's transversality (\ref{Griff}).

Finally a {\em polarization} of $(L,F)$ is a pairing of local systems $S: L\otimes_{\Q_Y} L
\to \Q_y$, that induces a polarization of Hodge structures on each stalk $L_y$ ($y\in Y$).
\end{defn}

For example in the geometric case above, one can get such a polarization on $L=R^nf_*\Q_X$
for $f: X\to Y$ a {\em projective smooth} morphism of complex algebraic (or analytic)
varieties. The existence of a polarization is needed for example for the following important result of Schmid (\cite{Sch}[thm.7.22]):

\begin{thm}[Rigidity]\label{rigid}
Let $Y$ be a connected complex manifold Zarisky open in a compact complex analytic manifold 
$\bar{Y}$, with $(L,F)$ a polarizable variation of pure Hodge structures on $Y$.
Then $H^0(Y,L)$ gets an induced Hodge structure such that the evaluation map
$H^0(Y,L)\to L_y$ is an isomorphism of Hodge structures for all $y\in Y$.
In particular the variation  $(L,F)$ is constant, if the underlying local system $L$ is constant.
\end{thm}

\subsection{Variation of mixed Hodge structures}
If one considers a morphism $f: X\to Y$ of complex algebraic varieties with $Y$ smooth, which is a topological fibration with possible singular or non-compact fiber, then the locally constant direct image sheaves $L=L^n:=R^nf_*\Q_X$ ($n\geq 0$) are {\em variations of mixed Hodge structures} in the sense of the following definitions.

\begin{defn}
Let $Y$ be a reduced complex analytic space. A {\em holomorphic family of mixed Hodge structures} on $Y$ consists of the following data:
\begin{enumerate}
\item a local system $L$ of rational vector spaces on $Y$ with finite dimensional stalks,
\item a finite decreasing {\em Hodge filtration} $F$ of $\LL=L\otimes_{\Q_Y} \OO_Y$ by holomorphic subbundles $F^p\LL$,
\item a finite increasing {\em weight filtration} $W$ of $L$ by local subsystems $W_nL$,
\end{enumerate}
such that the induced filtrations on $\LL_y\simeq L_y\otimes_{\Q}\C$ and $L_y$ define a
mixed Hodge structure on all stalks $L_y$ ($y\in Y$).

If $Y$ is a smooth complex manifold, then such a holomorphic family $(L,F,W)$ is called
a {\em variation of mixed Hodge structures}, if one has in addition
for the induced connection $\nabla: \LL\to \LL\otimes_{\OO_Y}\Omega_Y^1$ the Griffith's transversality (\ref{Griff}).

Finally $(L,F,W)$ is called {\em graded polarizable}, if the induced family (or variation)
of pure Hodge structures $Gr^W_nL$ (with the induced Hodge filtration $F$) is polarizable for all $n$. 
\end{defn}

With the obvious notion of morphisms, the categories $FmHs^{(p)}(Y)$ (or $VmHs^{(p)}(Y)$) of (graded polarizable) families (or variations) of mixed Hodge structures on $Y$ become abelian categories with
a tensor product $\otimes$ and duality $(\cdot)^{\vee}$. Again any such morphism is strictly compatible with the Hodge and weight filtrations. Moreover, one has 
for a holomorphic map $f: X\to Y$ (of complex manifolds) a functorial pullback
$$f^*: FmHs^{(p)}(Y)\to FmHs^{(p)}(X) \quad \text{(or $f^*: VmHs^{(p)}(Y)\to VmHs^{(p)}(X)$),}$$
comuting with tensor product $\otimes$ and duality $(\cdot)^{\vee}$.
On a point space $pt$ one gets just back  the category 
$$FmHs^{(p)}(pt)=VmHs^{(p)}(pt) =mHs^{(p)}$$ 
of (graded polarizable) mixed Hodge structures. 
Using the pullback under the constant map $k: Y\to pt$, we get the constant
family (or variation) of Tate Hodge structures $\Q_Y(n):=k^*\Q(n)$ on $Y$.

\subsection{Cohomological characteristic classes}
The Grothendieck group  $K^0_{an}(Y)$ of holomorphic vector bundles on the complex variety $Y$
is a commutative ring with multiplication induced by $\otimes$ and has a duality involution induced by $(\cdot)^{\vee}$. For a holomorphic map $f: X\to Y$ one has a
functorial pullback $f^*$ of rings with involutions. Similarly for $K^0_{an}(Y)[y^{\pm 1}]$,
if we extend the duality involution by 
$$([V]\cdot y^k)^{\vee}:=[V^{\vee}]\cdot (1/y)^k \:.$$

For a family (or variation) of mixed Hodge structures $(L,F,W)$ on $Y$ let us introduce 
the characteristic class 
\begin{equation}\label{MHC-coh}
MHC^y((L,F,W)):=\sum_p\; [Gr^p_F(\LL)]\cdot (-y)^p \in K^0_{an}(Y)[y^{\pm 1}]\:.
\end{equation}
Since morphisms of families (or variations) of mixed Hodge structures are strictly
compatible with the Hodge filtrations, we get an induced group homomorphism
of Grothendieck groups:
$$MHC^y: K_0(FmHs^{(p)}(Y))\to K^0_{an}(Y)[y^{\pm 1}] \quad 
\text{or $MHC^y: K_0(VmHs^{(p)}(Y))\to K^0_{an}(Y)[y^{\pm 1}]$.}$$

Note that  $MHC^{-1}((L,F,W))=[\LL]\in K^0_{an}(Y)$ is just the class of the associated
holomorphic vector bundle.
And for $Y=pt$ a point, we get back the $\chi_y$-genus:
$$\chi_y=MHC^y: K_0(mHs^{(p)})= K_0(FmHs^{(p)}(pt))\to K^0_{an}(pt)[y^{\pm 1}]
=\Z[y^{\pm 1}]\:.$$

\begin{thm}
The transformation 
$$MHC^y: K_0(FmHs^{(p)}(Y))\to K^0_{an}(Y)[y^{\pm 1}] \quad 
\text{or $MHC^y: K_0(VmHs^{(p)}(Y))\to K^0_{an}(Y)[y^{\pm 1}]$}$$
is contravariant functorial. It is a transformation of commutative rings with unit, i.e.
it commutes with products and respects the units: $MHC^y([\Q_Y(0)])=[\OO_Y]$.
Similarly it respects the duality involutions:
$$MHC^y([(L,F,W)^{\vee}])=\sum_p\; [(Gr^{-p}_F(\LL))^{\vee}]\cdot (-y)^p=
\left(MHC^y([(L,F,W)])\right)^{\vee}\:.$$
\end{thm}

\begin{example}\label{ex:smooth}
Let $f: X\to Y$ be a {\em proper smooth} morphism of complex algebraic varieties
or a {\em projective smooth} morphism of complex analytic varieties, so that
the higher direct image sheaf $L^n:=R^nf_*\Q_X$ ($n\geq 0$) with the induced Hodge filtration
as in (\ref{F-relDR}) defines a holomorphic family of pure Hodge structures on $Y$.
If $m$ is the complex dimension of the fibers, then $L_n=0$ for $n>2m$ so that one can define
$$[Rf_*\Q_X]:=\sum_{n=0}^{2m}\; (-1)^n\cdot [(R^nf_*\Q_X,F)] \in K_0(FmHs(Y))\:.$$
Then one gets by (\ref{F-relDR}):
\begin{equation}\begin{split} MHC^y([Rf_*\Q_X])
&= \sum_{p,q\geq 0}\;(-1)^{p+q}\cdot [R^qf_*\Omega^p_{X/Y}]\cdot (-y)^p\\
&= \sum_{p\geq 0}\; f_*[\Omega^p_{X/Y}]\cdot y^p\\
&=: f_*\left(\lambda_y(T^*_{X/Y})\right) \in K^0_{an}(Y)[y] \:.
\end{split}\end{equation}
Assume moreover that
\begin{enumerate}
\item[(a)] $Y$ is a connected complex manifold Zarisky open in a compact complex analytic manifold $\bar{Y}$,
\item[(b)] All direct images sheaves $L^n:=R^nf_*\Q_X$ ($n\geq 0$) are {\em constant}.
\end{enumerate}
Then one gets by the {\em rigidity} theorem \ref{rigid} (for $z\in Y$):
$$f_*\left(\lambda_y(T^*_{X/Y})\right)=\chi_y(\{f=z\})\cdot [\OO_Y]\in K^0_{an}(Y)[y]\:.$$
\end{example}

\begin{cor}[Multiplicativity]
Let $f: X\to Y$ be a {\em smooth} morphism of compact complex algebraic manifolds,
with $Y$ connected. Let $T^*_{X/Y}$ be the relative holomorphic cotangent bundle of the fibers,
fitting into the short exact sequence
$$0\to f^*T^*Y \to T^*X \to T^*_{X/Y}\to 0 \:.$$
Assume all direct images sheaves $L^n:=R^nf_*\Q_X$ ($n\geq 0$) are {\em constant}, i.e. $\pi_1(Y)$ acts trivially on the cohomology $H^*(\{f=z\})$ of the fiber.
Then one gets the multiplicativity of the $\chi_y$-genus (with $k: Y\to pt$ the constant map):
\begin{equation}\begin{split}\label{multipl-chiy}
\chi_y(X)&= (k\circ f)_* [\lambda_y(T^*X)]\\
&= k_*f_*\left( [\lambda_y(T^*_{X/Y})]\otimes f^*[\lambda_y(T^*Y)]\right)\\
&=k_*\left( \chi_y(\{f=z\})\cdot [\lambda_y(T^*Y)]\right)\\
&=\chi_y(\{f=z\})\cdot \chi_y(Y) \:.
\end{split}\end{equation}
\end{cor}

\begin{remark}
 The multiplicativity relation (\ref{multipl-chiy}) specializes for $y=1$ to the classical multiplicativity formula
$$sign(X)= sign(\{f=z\})\cdot sign(Y)$$
of Chern-Hirzebruch-Serre \cite{CHS} for the signature of an oriented fibration of smooth coherently oriented compact manifolds,
if $\pi_1(Y)$ acts trivially on the cohomology $H^*(\{f=z\})$ of the fiber. So it is a Hodge theoretic counterpart of this.
Moreover, the corresponding Euler characteristic formula for $y=-1$
$$\chi(X)= \chi(\{f=z\})\cdot \chi(Y)$$
is even true {\em without}  $\pi_1(Y)$ acting trivially on the cohomology $H^*(\{f=z\})$ of the fiber!

The Chern-Hirzebruch-Serre signature formula was motivational for many subsequent works which studied monodromy contributions
to invariants (genera and characteristic classes), e.g. see \cite{At, BCS, CMS, CMS0, CMS, CLMS, CLMS2, CS, MSc, Mey}. 
\end{remark}

Instead of working with holomorphic vector bundles, we can of course also use only the underlying topological complex vector bundles, which gives the forgetful transformation
$$For: K^0_{an}(Y)\to K^0_{top}(Y)\:.$$
Here the target can also be viewed as the even part of $\Z_2$-graded topological complex
K-cohomology. Of course, $For$ is contravariant functorial and commutes with
product $\otimes$ and duality $(\cdot)^{\vee}$. This duality induces a $\Z_2$-grading
on $K^0_{top}(Y)[1/2]$ by splitting into the (anti-)invariant part, and similarly for
$K^0_{an}(Y)[1/2]$.  Then the (anti-)invariant part of $K^0_{top}(Y)[1/2]$ can be identified 
with the even part of $\Z_4$-graded topological real K-theory $KO^0_{top}(Y)[1/2]$
(and $KO^2_{top}(Y)[1/2]$).\\

Assume now that $(L,F)$ is a holomorphic family of pure Hodge structures of weight $n$
on the complex variety $Y$, with a polarization $S: L\otimes_{\Q_Y} L\to \Q_Y$.
This induces an isomorphism of families of pure 
Hodge structures (of weight $n$):
$$L \simeq L^{\vee}(-n):=L^{\vee}\otimes \Q_Y(-n)\:.$$
So if we choose the isomorphism of rational local systems $\Q_Y(-n)=(2\pi i)^n\cdot \Q_Y\simeq \Q_Y$,
then the polarisation induces a $(-1)^n$-symmetric duality isomorphism $L\simeq L^{\vee}$
of the underlying local systems. And for such an (anti)symmetric selfdual local system $L$
Meyer \cite{Mey} has introduced a $KO$-characteristic class 
$$[L]_{KO}\in KO^0_{top}(Y)[1/2]\oplus KO^2_{top}(Y)[1/2])=K^0_{top}(Y)[1/2] $$
so that for $Y$ a compact oriented manifold of even real dimension $2m$ the following {\em twisted signature formula} is true:
\begin{equation}\label{tw-sig}
sign(H^m(Y,L))= \langle ch^*(\Psi^2([L]_{KO})), L^*(TM)\cap [M]\rangle \:.
\end{equation}
Here $H^m(Y,L)$ gets an induced (anti)symmetric duality, with $sign(H^m(Y,L)):=0$
in case of an antisymmetric pairing. Moreover $ch^*$ is the Chern character, $\Psi^2$
the second Adams operation and $L^*$ is the Hirzebruch-Thom $L$-class.\\

We now explain that $[L]_{KO}$ agrees up to some universal signs with
$For(MHC^1((L,F))$.
The underlying topological complex vector bundle of
$\LL$ has a natural real structure so that as a topological complex vector bundle one gets an orthogonal decomposition 
$$\LL =\oplus_{p+q=n}\; \HC^{p,q} \quad \text{with $\HC^{p,q}=F^p\LL\cap\overline{F^q\LL}
=\overline{\HC^{q,p}}$,}$$
with
\begin{equation}\label{MHC1}
For(MHC^1((L,F))= \sum_{p\;even,q}\; [\HC^{p,q}]- \sum_{p\;odd,q}\; [\HC^{p,q}]\:.
\end{equation}
If $n$ is even, then both sums of the right hand side in (\ref{MHC1}) are invariant under conjugation.
And $(-1)^{-n/2}\cdot S$ is by (\ref{positive}) positive resp. negative definite on the corresponding real vector bundle $(\oplus_{p\;even,q}\; \HC^{p,q})_{\R}$ resp.
$(\oplus_{p\;odd,q}\; \HC^{p,q})_{\R}$. So if we choose the pairing $(-1)^{n/2}\cdot S$
for the isomorphism $L\simeq L^{\vee}$, then this agrees with the splitting introduced by
Meyer \cite{Mey} in the definition of his $KO$-characteristic class $[L]_{KO}$ associated to this {\em symmetric} duality isomorphism of $L$:
$$For(MHC^1((L,F))= [L]_{KO} \in KO^0_{top}(Y)[1/2]\:.$$
Similarly, if $n$ is odd, both sums of the right hand side in (\ref{MHC1})  are exchanged under conjugation. If we choose the pairing $(-1)^{(n+1)/2}\cdot S$
for the isomorphism $L\simeq L^{\vee}$, then this agrees 
by definition \ref{pol} with the splitting introduced by
Meyer \cite{Mey} in the definition of his $KO$-characteristic class $[L]_{KO}$ associated to this {\em antisymmetric} duality isomorphism of $L$:
$$For(MHC^1((L,F))= [L]_{KO} \in KO^2_{top}(Y)[1/2]\:.$$

\begin{cor}
Let $(L,F)$ be a holomorphic family of pure Hodge structures of weight $n$
on the complex variety $Y$, with a polarization $S$ chosen.
Then the class $[L]_{KO}$ introduced in \cite{Mey} 
for the duality isomorphism coming from the pairing $(-1)^{n(n+1)/2}\cdot S$
is equal to 
$$For(MHC^1((L,F))=  [L]_{KO}\in KO^0_{top}(Y)[1/2]\oplus KO^2_{top}(Y)[1/2]
= K^0_{top}(Y)[1/2]\:.$$ 
It is therefore independent of the choice of the polarisation $S$. Moreover, this identification is functorial
under pullback and compatible with products (as defined in \cite{Mey}[p.26] for (anti)symmetic selfdual local systems).
\end{cor}

There are Hodge theoretic counterparts of the twisted signature formula (\ref{tw-sig}).
Here we formulate a corresponding K-theoretical result.
Let $(L,F,W)$ be a variation of mixed Hodge structures on the $m$-dimensional complex manifold $M$. Then 
$$H^n(M,L)\simeq H^n(M,DR(\LL))$$ 
gets an induced (decreasing) $F$ filtration
coming from the filtration of the holomorphic de Rham complex of the vector bundle
$\LL$ with its integrable connection $\nabla$:
$$ \begin{CD} DR(\LL)=[\LL @> \nabla >> \cdots @> \nabla >> \LL\otimes_{\OO_M} \Omega^m_M]
\end{CD}$$ 
(with $\LL$ in degree zero), defined by
\begin{equation}\label{dr-fil-ext} \begin{CD}
F^p DR(\LL) =[F^p\LL @> \nabla >> \cdots  @> \nabla >> F^{p-m}\LL\otimes_{\OO_M} \Omega^m_M] \:.
\end{CD}
\end{equation}
Note that here we are using the Griffith's transversality (\ref{Griff})!\\

The following result is due to Deligne and Zucker (\cite{Zu}[thm.2.9, lem.2.11]) in the case of a compact K\"{a}hler manifold,
whereas the case of a compact complex algebraic manifold follows from Saito's general results as explained in the next section.

\begin{thm}
Assume $M$ is a compact K\"{a}hler manifold or a compact complex algebraic manifold,
with  $(L,F,W)$ a graded polarizable variation of mixed (or pure) Hodge structures on $M$. Then $H^n(M,L)\simeq H^n(M,DR(\LL))$ gets an induced mixed (or pure) Hodge structure with
$F$ the Hodge filtration.
Moreover, the corresponding {\em Hodge to de Rham spectral-sequence}
degenerates at $E_1$ so that 
$$Gr_F^p(H^n(M,L)) \simeq H^{n}(M,Gr_F^pDR(\LL)) \quad \text{for all $n,p$.}$$
\end{thm}
Therefore one gets as a corollary (compare \cite{CLMS, CLMS2, MSc}):
\begin{equation}\begin{split}
\chi_y(H^*(M,L)) 
&= \sum_{n,p}\; (-1)^n\cdot dim_{\C}\left(H^{n}(M,Gr_F^pDR(\LL))\right) \cdot (-y)^p\\
&= \sum_p\; \chi\left(H^{*}(M,Gr_F^pDR(\LL))\right) \cdot (-y)^p\\
&= \sum_{p,i}\; (-1)^i\cdot \chi\left(H^{*}(M,Gr_F^{p-i}(\LL)\otimes_{\OO_M}\Omega^i_M)
\right)\cdot (-y)^p \\
&=k_*\left( MHC^y(L)\otimes  \lambda_y(T^*M)\right) \\
&=:\langle MHC^y(L), \lambda_y(T^*M)\cap [\OO_M] \rangle \in \Z[y^{\pm 1}]\:.
\end{split}
\end{equation}

\subsection{Good variation of mixed Hodge structures}
For later use let us introduce the following

\begin{defn}[good variation]\label{good}
Let $M$ be a complex algebraic manifold. A graded polarizable
variation of mixed Hodge structures
$(L,F,W)$ on $M$ is called good, if it is {\em admissible} in the sense of Steenbrink-Zucker \cite{SZ} and Kashiwara \cite{Ka},
with {\em quasi-unipotent monodromy} at infinity, i.e. with respect to a compactification $\bar{M}$ of $M$ by a compact complex algebraic manifold $\bar{M}$, with complement $D:=\bar{M}\backslash M$ a normal crossing divisor with smooth irreducible components.
\end{defn}

\begin{example}[pure and geometric variations] \label{good-var}
Two important examples for such a good variation of mixed Hodge structures are the following:
\begin{enumerate}
 \item A polarizable variation of {\em pure} Hodge structures is always admissible
  by a deep theorem of Schmid \cite{Sch}[thm.6.16]. So it is good precisely when it has quasi-unipotent monodromy at infinity.
\item Consider a morphism $f: X\to Y$ of complex algebraic varieties with $Y$ smooth, which is a topological fibration with possible singular or non-compact fiber. Then the locally constant direct image sheaves $R^nf_*\Q_X$ and $R^nf_!\Q_X$ ($n\geq 0$) are good variations of mixed Hodge structures (compare with remark \ref{geometric}).
\end{enumerate}
\end{example}

This class of good variations on $M$ is again an abelian category $VmHs^g(M)$
stable under tensor product $\otimes$,
duality $(\cdot)^{\vee}$ and pullback $f^*$ for $f$ an algebraic morphism of complex algebraic manifolds. Moreover, in this case all vector bundles $F^p\LL$ of the Hodge filtration carry the structure of a unique underlying complex algebraic vector bundle (in the Zariski topology),
so that the characteristic class transformation $MHC^y$ can be seen as a natural contravariant transformation of rings with involution
$$MHC^y: K_0(VmHs^g(M)) \to K^0_{alg}(M)[y^{\pm 1}] \:.$$

In fact, consider a (partial) compactification $\bar{M}$ of $M$ as above, with
$D:=\bar{M}\backslash M$ a normal crossing divisor with smooth irreducible components
and $j: M\to \bar{M}$ the open inclusion.
Then the holomorphic vector bundle $\LL$ with integrable connection $\nabla$ corresponding to
$L$ has a unique {\em canonical Deligne extension} $(\overline{\LL},\overline{\nabla})$ to a holomorphic
vector bundle $\overline{\LL}$ on $\bar{M}$, with {\em meromorphic} integrable connection 
\begin{equation}\label{D-ext}
\overline{\nabla}: \overline{\LL}\to 
\overline{\LL}\otimes_{\OO_{\bar{M}}} \Omega^1_{\bar{M}}(log(D))
\end{equation}
having {\em logarithmic poles} along $D$. Here the {\em residues} of $\overline{\nabla}$
along $D$ have real eigenvalues, since $L$ has {\em quasi-unipotent monodromy} along $D$.
And the canonical extension is characterized by the property, that all these
eigenvalues are in the half-open intervall $[0,1[\;$.
Moreover, also the Hodge filtration $F$ of $\LL$ extends uniquely to a filtration
$\bar{F}$ of $\overline{\LL}$ by holomorphic subvector bundles
$$F^p\overline{\LL}:=j_*(F^p\LL) \cap \overline{\LL} \subset j_*\LL \:,$$ 
since $L$ is {\em admissible} along $D$. Finally the Griffith's transversality extends to
\begin{equation}\label{Griff-alg}
\overline{\nabla}(F^p\overline{\LL})\subset F^{p-1}\overline{\LL}\otimes_{\OO_{\bar{M}}} \Omega^1_{\bar{M}}(log(D)) \quad\text{for all $p$.}
\end{equation}
For more details see \cite{De}[prop.5.4] and \cite{PS}[sec.11.1, sec.14.4].\\

If we choose $\bar{M}$ as a compact algebraic manifold, then we can apply Serre's GAGA theorem 
to conclude that $\overline{\LL}$ and all $F^p\overline{\LL}$ are {\em algebraic} vector bundles, with $\overline{\nabla}$ an {\em algebraic} meromorphic connection.

\begin{remark} \label{rem-del}
The canonical Deligne extension $\overline{\LL}$ (as above) with its Hodge filtration $F$ has the following compabilities
(compare \cite{De}[part II]):
\begin{description}
\item[smooth pullback] Let $f:\bar{M'} \to \bar{M}$ be a smooth morphism so that $D':=f^{-1}(D)$
is also a normal crossing divisor with smooth irreducible components on $\bar{M'}$ with complement $M'$. Then one has 
\begin{equation}
f^*\left(\overline{\LL}\right) \simeq \overline{f^*\LL} \quad \text{and}
 f^*\left(F^{p}\overline{\LL}\right) \simeq F^p\overline{f^*\LL}\:\:\text{for all $p$.}
\end{equation}
 \item[exterior product]  
Let $L$ and $L'$ be  two good variations on $M$ and $M'$. Then their canonical Deligne extensions satisfy
$$\overline{\LL\boxtimes_{\OO_{M\times M'}} \LL'}
\simeq \overline{\LL}\boxtimes_{\OO_{\bar{M}\times \bar{M}'}} \overline{\LL'}\:,$$
since the residues of the corresponding  meromorphic connections are compatible. 
Then one has for all $p$:
\begin{equation}
F^p\left(\overline{\LL\boxtimes_{\OO_{M\times M'}} \LL'}\right)\simeq \oplus_{i+k=p}\; \left(F^i\overline{\LL}\right)\boxtimes_{\OO_{\bar{M}\times \bar{M}'}} \left(F^k\overline{\LL'}\right) \:.
\end{equation}
\item[tensor product] In general the canonical Deligne extensions of two 
good variations $L$ and $L'$ on $M$ are {\em not} compatible with tensor products,
because of the choice of different residues for the corresponding  meromorphic connections. This problem doesn't appear if one of these variations,
lets say $L'$, is already defined on $\bar{M}$. 
Let $L$ resp. $L'$ be a good variation on $M$ resp. $\bar{M}$. Then their canonical Deligne extensions satisfy
$$\overline{\LL\otimes_{\OO_M} (\LL'|M)}
\simeq \overline{\LL}\otimes_{\OO_{\bar{M}}} \LL'\:,$$
and one has for all $p$:
\begin{equation}
F^p\left(\overline{\LL\otimes_{\OO_M} (\LL'|M)}\right)\simeq \oplus_{i+k=p}\; \left(F^i\overline{\LL}\right)\otimes_{\OO_{\bar{M}}} \left(F^k\LL'\right) \:.
\end{equation}
\end{description} 
\end{remark}

Let $\bar{M}$ be a partial compactification of $M$ as before, i.e. we don't assume that $\bar{M}$ is compact, with $m:=dim_{\C}(M)$.
Then the {\em logarithmic de Rham complex}
$$ \begin{CD} DR_{log}\left(\overline{\LL}\right)
:=[\overline{\LL} @> \overline{\nabla} >> \cdots @> \overline{\nabla} >> \overline{\LL}\otimes_{\OO_{\bar{M}}} \Omega^m_{\bar{M}}(log(D))]
\end{CD}$$ 
(with $\overline{\LL}$ in degree zero) is quasi-isomorphic to $Rj_*L$, so that
$$H^*(M,L) \simeq H^*\left(\bar{M},DR_{log}\left(\overline{\LL}\right)\right) \:.$$
So these cohomology groups get an induced (decreasing) $F$-filtration coming from the filtration
\begin{equation}\label{drlog-fil} \begin{CD}
F^p DR_{log}\left(\overline{\LL}\right) =
[F^p\overline{\LL} @> \overline{\nabla} >> \cdots @> \overline{\nabla} >> F^{p-m}\overline{\LL}\otimes_{\OO_{\bar{M}}} \Omega^m_{\bar{M}}(log(D))] \:.
\end{CD}
\end{equation}

For $\bar{M}$ a compact algebraic manifold, this is again the Hodge filtration of an induced mixed Hodge structure on $H^*(M,L)$
(compare with corollary \ref{MHM3}).

\begin{thm}\label{j*VMHS}
Assume $\bar{M}$ is a smooth algebraic compactification of the algebraic manifold $M$ with the complement $D$ a normal crossing divisor with smooth irreducible components.
Let $(L,F,W)$ be a good variation of mixed Hodge structures on $M$. Then $H^n(M,L)\simeq H^*\left(\bar{M},DR_{log}\left(\overline{\LL}\right)\right)$ gets an induced mixed Hodge structure with $F$ the Hodge filtration.
Moreover, the corresponding {\em Hodge to de Rham spectral-sequence}
degenerates at $E_1$ so that 
$$Gr_F^p(H^n(M,L)) \simeq H^{n}\left(M,Gr_F^pDR_{log}\left(\overline{\LL}\right)\right) \quad \text{for all $n,p$.}$$
\end{thm}
Therefore one gets as a corollary (compare \cite{CLMS, CLMS2, MSc}):
\begin{equation}\label{pairing1}\begin{split}
\chi_y(H^*(M,L)) 
&= \sum_{n,p}\; (-1)^n\cdot dim_{\C}\left(H^{n}\left(M,Gr_F^pDR_{log}\left(\overline{\LL}\right)\right)\right) \cdot (-y)^p\\
&= \sum_p\; \chi\left(H^{*}\left(M,Gr_F^pDR_{log}\left(\overline{\LL}\right)\right)\right) \cdot (-y)^p\\
&= \sum_{p,i}\; (-1)^i\cdot \chi\left(H^{*}\left(M,Gr_F^{p-i}\left(\overline{\LL}\right)\otimes_{\OO_{\bar{M}}}\Omega^i_{\bar{M}}(log(D)) \right)
\right)\cdot (-y)^p \\
&=:\langle MHC^y(Rj_*L), \lambda_y\left(\Omega^1_{\bar{M}}(log(D))  \right)\cap [\OO_{\bar{M}}] \rangle \in \Z[y^{\pm 1}]\:.
\end{split}
\end{equation}

Here we use the notion
\begin{equation}\label{MHCj*-coh}
MHC^y(Rj_*L):=\sum_p\; [Gr^p_F\left(\overline{\LL}\right)]\cdot (-y)^p \in K^0_{alg}(\bar{M})[y^{\pm 1}]\:.
\end{equation}

Remark \ref{rem-del} then implies the

\begin{cor}\label{Cor-j*}
 Let $\bar{M}$ be a smooth algebraic partial compactifiction of the algebraic manifold $M$ with the complement $D$ a normal crossing divisor with smooth irreducible components.
Then $MHC^y(Rj_*(\cdot))$ induces a transformation 
$$MHC^y(j_*(\cdot)): K_0(VmHs^{g}(M))\to K^0_{alg}(\bar{M})[y^{\pm 1}] \:.$$
\begin{enumerate}
\item This is contravariant functorial for a smooth morphism $f: \bar{M}'\to \bar{M}$ of such partial compactifications, i.e. 
$$f^*\left(MHC^y(j_*(\cdot))\right) \simeq MHC^y\left(j'_*(f^*(\cdot))\right) \:.$$
\item 
It commutes with exterior products for two good variations $L,L'$:
$$MHC^y\left((j\times j')_*[(L\boxtimes_{\Q_{M\times M'}} L']\right)= 
MHC^y(j_*[L]) \boxtimes MHC^y(j'_*[(L']) \:.$$
\item
Let $L$ resp. $L'$ be a good variation on $M$ resp. $\bar{M}$. Then $MHC^y(j_*[\cdot])$ is multiplicative in the sense that
$$MHC^y\left(j_*[(L\otimes_{\Q_M} (L'|M)]\right)= 
MHC^y(j_*[L]) \otimes MHC^y([L']) \:.$$
\end{enumerate}
\end{cor}

\section{Calculus of mixed Hodge modules}
\subsection{Mixed Hodge modules.} Before discussing extensions of
the characteristic cohomology classes $MHC^y$  to the singular setting, we need to briefly
recall some aspects of Saito's theory \cite{Sa0,Sa1,Sa2, Sa4, Sa5} of algebraic mixed Hodge
modules, which play the role of singular extensions of good variations of mixed Hodge structures.\\

To each complex algebraic variety $Z$, Saito associated a
category $MHM(Z)$ of {\em algebraic mixed Hodge modules} on $Z$ (cf.
\cite{Sa0,Sa1}). If $Z$ is smooth, an object of this category
consists of an algebraic (regular) holonomic $D$-module $(\M,F)$ with a good filtration $F$
together with a perverse sheaf $K$ of rational vector spaces, both endowed a finite
increasing filtration $W$ such that
$$\alpha: DR(\M)^{an}\simeq K\otimes_{\Q_Z} \C_Z \quad \text{is compatible with $W$}$$
under the Riemann-Hilbert correspondence coming from the (shifted) analytic  de Rham complex
(with $\alpha$ a chosen isomorphism). Here we use left $D$-modules, and the sheaf $\DC_Z$ of algebraic differential operators on $Z$ has the increasing filtration $F$ with $F_i\DC_Z$ given by the differential operators of degree
$\leq i$ ($i\in \Z$). Then a {\em good} filtration $F$ of the algebraic holonomic $D$-module $\M$ is given by a 
bounded from below, increasing and exhaustive filtration $F_p\M$ by {\em coherent} algebraic $\OO_Z$-modules such that
\begin{equation}\label{D-filt}
F_i\DC_Z\left( F_p\M\right) \subset F_{p+i}\M\; \quad \text{for all $i,p$, and this is an equality for $i$ big enough.}
\end{equation}

In general, for a singular variety $Z$ one works
with suitable local embeddings into manifolds and corresponding
filtered $D$-modules supported on $Z$. In addition, these objects
are required to satisfy a long list of complicated properties (not needed here).
The {\em forgetful} functor $rat$ is defined as
$$rat:  MHM(Z)\to Perv(\Q_Z);\: \left((\M,F),K,W\right)\mapsto K\:.$$

\begin{thm}[M. Saito]\label{MHM1}
$MHM(Z)$ is an abelian category with $rat:  MHM(Z)\to Perv(\Q_Z)$
exact and {\em faithful}. It extends to a functor 
$$rat: D^bMHM(Z) \to D^b_c(\Q_Z)$$ 
to the derived category of complexes of $\Q$-sheaves with algebraically constructible
cohomology. There are functors 
$$f_*, \;f_! ,\; f^*,\; f^!,\; \otimes,\; \boxtimes,\; \DC \quad \text{on $D^bMHM(Z)$} \:, $$ 
which are
``lifts" via $rat$ of the similar (derived) functors defined on $D^b_c(\Q_Z)$, with $(f^*,f_*)$ and $(f_!,f^!)$ also pairs of {\em adjoint} functors.
One has a natural map $f_!\to f_*$, which is an isomorphism for $f$ proper.
Here $\DC$ is a duality involution $\DC^2\simeq id$ ``lifting" the Verdier duality functor, with
$$\DC\circ f^* \simeq f^!\circ \DC \quad \text{and} \quad
\DC\circ f_* \simeq f_!\circ \DC \:.$$
\end{thm}

Compare with \cite{Sa1}[thm.0.1 and sec.4] for more details (as well as with \cite{Sa4} for a more general formal abstraction).
The usual truncation $\tau_{\leq}$ on $D^bMHM(Z)$
corresponds to the {\em perverse truncation} ${^p\tau}_{\leq}$ on
$D^b_c(Z)$ so that
$$rat\circ H= \;^p\HC\circ rat\:,$$
where $H$ stands for the cohomological functor in $D^bMHM(Z)$ and $\;^p\HC$ denotes the perverse cohomology
(always with respect to the self-dual middle perversity).

\begin{example}\label{smooth}
Let $M$ be a complex algebraic manifold of pure complex dimension $m$,
with $(L,F,W)$ a good variation of mixed Hodge structures on $M$.
Then $\LL$ with its integrable connection $\nabla$ is a holonomic (left) $D$-module
with $\alpha: DR(\LL)^{an}\simeq L[m]$, where this time we use the shifted de Rham complex
$$\begin{CD}
DR(\LL):=[\LL @> \nabla >> \cdots @> \nabla >> \LL\otimes_{\OO_M} \Omega^m_M]
\end{CD}$$ 
with $\LL$ in degree $-m$, so that $DR(\LL)^{an}\simeq L[m]$ is a perverse sheaf on $M$.
The filtration $F$ induces by Griffith's transversality (\ref{Griff}) a good filtration
$F_p(\LL):=F^{-p}\LL$ as a filtered  $D$-module. 
As explained before, this comes from an underlying algebraic filtered $D$-module.
Finally $\alpha$ is compatible with the
induced filtration $W$ defined by 
$$W^i(L[m]):=W^{i-m}L[m] \quad \text{ and} \quad  
W^i(\LL):=(W^{i-m}L)\otimes_{\Q_M} \OO_M \:.$$
 And this defines a mixed Hodge module
$\M$ on $M$, with $rat(\M)[-m]$ a local system on $M$.  
\end{example}

A mixed Hodge module $\M$ on the pure $m$-dimensional complex algebraic manifold
$M$ is called {\em smooth}, if $rat(\M)[-m]$ is a local system on $M$. Then this example corresponds to
\cite{Sa1}[thm.0.2], whereas the next theorem corresponds to \cite{Sa1}[thm.3.27 and rem. on p.313]:

\begin{thm}[M. Saito]\label{MHM2}
Let $M$ be a pure $m$-dimensional complex algebraic manifold.
Associating to a good variation of mixed Hodge structures $\VB=(L,F,W)$ on $M$
the mixed Hodge module $\M:=\VB_H$ as in example (\ref{smooth}) defines an
equivalence of categories
$$MHM(M)_{sm}\simeq VmHs^g(M)$$
between the categories of smooth mixed Hodge modules $MHM(M)_{sm}$ and
good variation of mixed Hodge structures on $M$.
This commutes with exterior product $\boxtimes$ and pullback 
$$f^*: VmHs^g(M)\to VmHs^g(M') \quad \text{resp.}\quad f^*[m'-m]: MHM(M)\to MHM(M')$$
for an algebraic morphism of smooth algebraic manifolds $M,M'$ of dimension $m,m'$.
For $M=pt$ a point, one gets in particular an equivalence
$$MHM(pt)\simeq mHs^p \:.$$
\end{thm}

\begin{remark}\label{geometric}
These two theorems explain why a geometic variations of mixed Hodge structures as in
Example \ref{good-var}(2) is good.
\end{remark}

By the last identification of the theorem, there exists a unique Tate object  $\Q^H(n)
\in MHM(pt)$ such that $rat(\Q^H(n))=\Q(n)$ and $\Q^H(n)$ is of type $(-n,-n)$:
$$MHM(pt)\ni \Q^H(n) \simeq \Q(n)\in mHs^p \:.$$
 For a complex variety $Z$ with constant map $k: Z\to pt$, define 
$$\Q_Z^H(n):=k_Z^*\Q^H(n) \in D^bMHM(Z), \quad \text{with $rat(\Q_Z^H(n))=\Q_Z(n)$.}$$ 
So tensoring with $ \Q_Z^H(n)$ defines the Tate twist $\cdot (n)$  of mixed Hodge modules. 
To simplify the notations, let $\Q_Z^H:=\Q_Z^H(0)$.
If $Z$ is \emph{smooth} of complex dimension $n$ then
$\Q_Z[n]$ is perverse on $Z$, and $\Q_Z^H[n]\in MHM(Z)$ is a single
mixed Hodge module, explicitly described by
$$\Q_Z^H[n]=((\OO_Z, F), \Q_Z[n], W), \quad \text{with  $gr^F_i=0=gr^W_{i+n}$ for all $ i \neq 0$.}$$ 

It follows from the definition that every $\M \in MHM(Z)$ has a finite
increasing {\em weight filtration} $W$ so that the functor $M \to Gr^W_kM$
is exact. We say that 
$\M \in D^bMHM(Z)$ has {\em weights $\leq n$ (resp. $\geq n$)} if
$Gr_j^WH^iM=0$ for all $j>n+i$ (resp. $j<n+i$). $\M$ is called {\em pure of weight $n$},
if it has weights both $\leq n$ and $\geq n$. For the following results compare with
\cite{Sa1}[prop.2.26 and (4.5.2)]:

\begin{prop} 
 If $f$ is a map of algebraic varieties, then $f_!$ and
$f^*$ preserve weight $\leq n$, and $f_*$ and $f^!$ preserve weight
$\geq n$. If $f$ is smooth of pure complex fiber dimension $m$, then
$f^!\simeq f^*[2m](m)$ so that $f^*,f^!$ preserve pure objects for $f$ smooth.
Moreover, if $\M \in D^bMHM(X)$ is pure and $f:X \to
Y$ is proper, then $f_*\M \in D^bMHM(Y)$ is pure of the same weight
as $\M$. 

Similarly the duality functor $\DC$ exchanges ``weight $\leq n$'' and
 ``weight $\geq -n$'', in particular it preserves pure objects.
Finally let $j: U\to Z$ be the inclusion of a Zariski open subset. Then the {\em
intermediate extension} functor
\begin{equation}
j_{!*}: MHM(U) \to MHM(Z):\; \M\mapsto 
Im\left( H^0(j_!\M) \to H^0(j_*(\M)\right)
\end{equation} 
preserves  weight $\leq n$ and $\geq n$, in particular it preserves pure objects
(of weight $n$). 
\end{prop}

We say that $\M \in D^bMHM(Z)$ is supported on $S\subset Z$ if and only if
$rat(\M)$ is supported on $S$.
There are the abelian subcategories $MH(Z,k)^p \subset MHM(Z)$ of pure
Hodge modules of weight $k$, which in the algebraic context are assumed to be
polarizable (and extendable at infinity). 

For each $k \in \Z$, the
abelian category $MH(Z,k)^p$ is semi-simple, in the sense that every
pure Hodge module on $Z$ can be uniquely written as a finite direct
sum of pure Hodge modules with strict support in irreducible
closed subvarieties of $Z$. Let $MH_S(Z,k)^p$ denote the subcategory
of {\em pure Hodge modules of weight $k$ with strict support in
$S$}. Then every $\M \in MH_S(Z,k)^p$ is generically a good
variation of Hodge structures $\VB_U$ of weight $k-d$ ($d:=dim\;S$) on a Zariski dense smooth open subset
$U \subset S$ (i.e. $\VB_U$ is polarizable with quasi-unipotent monodromy at infinity).
This follows from theorem \ref{MHM2} and the fact, that a perverse sheaf is 
generically a shifted local system on a smooth dense Zariski open subset $U\subset S$.
Conversely, every such good variation of Hodge structures $\VB$ on such an $U$ 
corresponds by theorem \ref{MHM2} to a pure Hodge module $\VB_H$ on $U$, which can
be extended in an unique way to a pure Hodge module $j_{!*}\VB_H$ on $S$ with strict support. Under this
correspondence, for $M \in MH_S(Z,k)^p$ we have that
$$rat(\M)=IC_S(\VB)$$
is the {\em twisted intersection cohomology complex} for $\VB$ the corresponding variation of Hodge
structures. Similarly
\begin{equation}\label{duality-IC}
 \DC(j_{!*}\VB_H) \simeq j_{!*}(\VB^{\vee}_H) (d)\:.
\end{equation}

Moreover, a {\em polarization} of $\M \in MH_S(Z,k)^p$ corresponds to an isomorphism
of Hodge modules (compare \cite{PS}[def.14.35, rem.14.36])
\begin{equation}\label{pol-HM}
S: \M \simeq \DC(\M)(-k)\:,
\end{equation}
whose restriction to $U$ gives a polarization of $\VB$.
In particular it induces a self-duality  isomorphism
$$S: rat(\M) \simeq \DC(rat(\M))(-k) \simeq \DC(rat(\M))$$
of the underlying twisted intersection cohomology complex, if an isomorphism $\Q_U(-k)\simeq \Q_U$ is chosen.\\

So if $U$ is smooth of pure complex dimension $n$, then $\Q_U^H[n]$ is a pure 
Hodge module of weight $n$.  If moreover $j: U \hookrightarrow
Z$ is a Zariski-open dense subset in $Z$, then the {\em intermediate
extension } $j_{!*}$ for mixed Hodge modules (cf. also with \cite{BBD}) preserves the weights. This
shows that if $Z$ is a complex algebraic variety of pure dimension
$n$ and $j: U \hookrightarrow Z$ is the inclusion of a smooth
Zariski-open dense subset then the intersection cohomology module
$IC_Z^H:=j_{!*}(\Q_U^H[n])$ is pure of weight $n$, with underlying
perverse sheaf $rat(IC_Z^H)=IC_Z$.\\

Note that the stability of a pure object $\M\in MHM(X)$ under a proper morphism 
$f: X\to Y$ implies the famous {\em decomposition theorem} of \cite{BBD} in the
context of pure Hodge modules (\cite{Sa1}[(4.5.4) on p.324]):
\begin{equation}\label{decomp}
f_*\M \simeq \oplus_i\; H^if_*\M[-i]\:, \quad \text{with $H^if_*\M$ semi-simple for all $i$.}
\end{equation}

Assume $Y$ is pure-dimesional, with $f: X\to Y$ a {\em resolution of singularities},
i.e. $X$ is smooth with $f$ a proper morphism, which generically  is an isomorphism
on some Zariski dense open subset $U$.
Then $\Q^H_X$ is pure, since $X$ is smooth, and $IC^H_Y$ has to be the direct summand of $H^0f_*\Q^H_X$ which corresponds to $\Q^H_U$. 
\begin{cor}
Assume $Y$ is pure-dimesional, with $f: X\to Y$ a {\em resolution of singularities}.
Then $IC^H_Y$ is a direct summand of $f_*\Q^H_X \in D^bMHM(Y)$.
\end{cor}

Finally we get the following results about the existence of a mixed Hodge structure
on the cohomology (with compact support) $H^i_{(c)}(Z,M)$ for $\M\in D^bMHM(Z)$.

\begin{cor}\label{MHM3}
Let $Z$ be a complex algebraic variety with constant map $k: Z\to pt$.
Then the cohomology (with compact support) $H^i_{(c)}(Z,\M)$ of $\M\in D^bMHM(Z)$ gets an induced graded polarizable mixed Hodge structure:
$$H^i_{(c)}(Z,\M)=H^i(k_{*(!)}M)\in MHM(pt)\simeq mHs^p\:.$$
In particular:
\begin{enumerate}
\item The rational cohomology (with compact support) $H^i_{(c)}(Z,\Q)$ 
of $Z$ gets an induced graded polarizable mixed Hodge structure by:
$$H^i(Z,\Q)= rat(H^i(k_*k^*\Q^H)) \quad \text{and} \quad 
H^i_{c}(Z,\Q)= rat(H^i(k_!k^*\Q^H)) \:.$$
\item Let $\VB_U$ be a good variation of mixed Hodge structures on a
smooth pure $n$-dimensional complex variety $U$, which is Zariski open and dense in
a variety $Z$, with $j: U\to Z$ the open inclusion. Then the global twisted Intersection cohomology (with compact support)
$$IH_{(c)}^i(Z,\VB):=H^i_{(c)}(Z, IC_Z(\VB)[-n])$$ 
gets a mixed Hodge structure by
$$IH_{(c)}^i(Z,\VB)=H^i(k_{*(!)}IC_Z(\VB)[-n]) =H^i(k_{*(!)}j_{!*}(\VB)[-n]) \:.$$ 
If $Z$ is compact, with $\VB$ a polarizable variation
of pure Hodge structures of weight $w$, then also $IH^i(Z,\VB)$
has a (polarizable) pure Hodge structure of weight $w+i$.
\item Let $\VB$ be a good variation of mixed Hodge structures on a
smooth (pure dimensional) complex manifold $M$, which is Zariski open and dense in
complex algebraic manifold $\bar{M}$, with complement $D$ a normal crossing divisor
with smooth irreducible components. Then $H^i(M,\VB)$
gets a mixed Hodge structure by
$$H^i(M,\VB)\simeq H^i(\bar{M},j_*\VB)\simeq H^i(k_*j_*\VB)\:,$$
with $j: U\to Z$ the open inclusion.
\end{enumerate}
\end{cor}

\begin{remark} Let us point out some important properties of these mixed Hodge structures:
\begin{enumerate}
\item By a deep theorem of Saito (\cite{Sa5}[thm.0.2,cor.4.3]), the mixed Hodge structure on $H^i_{(c)}(Z,\Q)$ defined as above coincides with the classical mixed Hodge structure constructed by Deligne (\cite{De1, De3}).
\item Assume we are in the context of (3) above with $Z=\bar{M}$ projective
and $\VB$ a good variation of pure Hodge structures on $U=M$.
Then the pure Hodge structure of (2) on the global Intersection cohomology $IH^i(Z,\VB)$ agrees with that of  \cite{CKS, KK} 
defined in terms of $L^2$-cohomology with respect to a K\"{a}hler metric with Poincar\'{e}
singularities along $D$ (compare \cite{Sa1}[rem.3.15]). The case of a $1$-dimensional complex
algebraic curve $Z=\bar{M}$ due to Zucker \cite{Zu}[thm.7.12] is used in the work of Saito \cite{Sa0}[(5.3.8.2)]
in the proof of the stability of pure Hodge modules under projective morphisms \cite{Sa0}[thm.5.3.1]
(compare also with the detailed discussion of this $1$-dimensional case in \cite{Sab}).
\item Assume we are in the context of (3) above with $\bar{M}$ compact. Then the mixed Hodge structure on $H^i(M,\VB)$ is the one of theorem \ref{j*VMHS}, whose Hodge filtration $F$ 
comes from the filtered logarithmic de Rham complex (compare \cite{Sa1}[sec.3.10, prop.3.11]). 
\end{enumerate}
\end{remark}

\subsection{Grothendieck groups of algebraic mixed Hodge modules.}\label{Grot} 
In this section, we describe the functorial calculus of Grothendieck groups
of algebraic mixed Hodge modules. Let $Z$ be a complex algebraic
variety. By associating to (the class of) a complex the alternating
sum of (the classes of) its cohomology objects, we obtain the
following identification (e.g. compare [\cite{KS}, p. 77],
[\cite{Sc}, Lemma 3.3.1])
\begin{equation} K_0(D^bMHM(Z))=K_0(MHM(Z)).
\end{equation}
In particular, if $Z$ is a point, then
\begin{equation} K_0(D^bMHM(pt))=K_0(mHs^p),
\end{equation}
and the latter is a commutative ring with respect to the tensor
product, with unit $[\Q^H]$.  Then we have
for any complex $\M^{\bullet} \in D^bMHM(Z)$  the
identification
\begin{equation}\label{i1}
[\M^{\bullet}]=\sum_{i \in \Z} (-1)^i [H^i(\M^{\bullet})] \in
K_0(D^bMHM(Z)) \cong K_0(MHM(Z)).
\end{equation}
In particular, if for any $\M \in MHM(Z)$ and $k \in \Z$ we regard
$\M[-k]$ as a complex concentrated in degree $k$, then
\begin{equation}\label{i2}
\left[ \M[-k] \right]= (-1)^k [\M] \in K_0(MHM(Z)).
\end{equation}
All functors $f_*$, $f_!$, $f^*$, $f^!$, $\otimes$, $\boxtimes$, $\DC$
induce corresponding functors on $K_0(MHM(\cdot))$. Moreover,
$K_0(MHM(Z))$ becomes a $K_0(MHM(pt))$-module, with the
multiplication induced by the exact exterior product with a point space:
$$\boxtimes : MHM(Z) \times MHM(pt) \to MHM(Z \times \{pt\}) \simeq
MHM(Z).$$ Also note that 
$$\M \otimes \Q^H_Z \simeq \M \boxtimes
\Q^H_{pt} \simeq \M$$ 
for all $\M \in MHM(Z)$. Therefore,
$K_0(MHM(Z))$ is a unitary $K_0(MHM(pt))$-module. The functors
$f_*$, $f_!$, $f^*$, $f^!$ commute with exterior products (and $f^*$
also commutes with the tensor product $\otimes$), so that the
induced maps at the level of Grothendieck groups $K_0(MHM(\cdot))$
are $K_0(MHM(pt))$-linear. Similarly $\DC$ defines an involution on $K_0(MHM(\cdot))$.
Moreover, by the functor
$$rat:K_0(MHM(Z)) \to K_0(D^b_c(\Q_Z)) \simeq K_0(Perv(\Q_Z)),$$ all these transformations
lift the corresponding transformations from the (topological) level
of Grothendieck groups of constructible (or perverse) sheaves.

\begin{remark}
The Grothendieck group $K_0(MHM(Z))$ has two different type of generators:
\begin{enumerate}
\item It is generated by the classes of pure Hodge modules $[IC_S(\VB)]$ with strict
support in an irreducible complex algebraic subset $S\subset Z$, with $\VB$
a good variation of (pure) Hodge structures on a dense Zariski open smooth subset $U$  of $S$.
These generators behave well under duality.
\item It is generated by the classes $f_*[j_*\VB]$, with $f: \bar{M}\to Z$ a proper morphisms from the smooth complex algebraic manifold $\bar{M}$,
$j: M\to \bar{M}$ the inclusion of a Zariski open and dense subset $M$, with complement $D$ a normal crossing divisor with smooth irreducible components, and $\VB$ a good variation of
mixed (or if one wants also pure) Hodge structures on $M$. These generators  will
be used in the next section about characteristic classes of mixed Hodge modules.
\end{enumerate}
\end{remark}

Here (1) follows from the fact, that a mixed Hodge module has a finite weight filtration,
whose graded pieces are pure Hodge modules, i.e. are finite direct sums of pure Hodge modules
$IC_S(\VB)$ with strict support $S$ as above. (2) follows by induction from resolution of singularities and
from the existence of a standard distinguished triangle associated to a closed inclusion.\\

Let $i: Y\to Z$ be a closed inclusion of complex algebraic varieties
 with open complement $j: U=Z\backslash Y\to Z$. Then one has by Saito's
work \cite{Sa1}[(4.4.1)]   the following functorial distinguished triangle in $D^bMHM(Z)$: 
\begin{equation}\label{triangle}
 \begin{CD}
  j_!j^* @> ad_j >> id @> ad_i >> i_*i^* @> [1]>> \:.
 \end{CD}
\end{equation}
Here the maps $ad$ are the adjunction maps, with $i_*=i_!$ since $i$ is proper.
If $f: Z \to X$ is a complex algebraic morphism, then we can apply $f_!$ to get another 
distinguished triangle
\begin{equation}\label{triangle2}
 \begin{CD}
  f_!j_!j^*\Q^H_Z @> ad_j >> f_!\Q^H_Z @> ad_i >> f_!i_!i^*\Q^H_Z @> [1]>> \:.
 \end{CD}
\end{equation}
On the level of Grothendieck groups, we get the important {\em additivity relation}
\begin{equation}
f_![\Q^H_Z] = (f\circ j)_![\Q^H_U] + (f\circ i)_![\Q^H_Y] \in K_0(D^bMHM(X))=K_0(MHM(X))\:.
\end{equation}
 
\begin{cor}\label{chiHdg}
One has a natural group homomorphism
$$\chi_{Hdg}: K_0(var/X) \to K_0(MHM(X)); [f:Z\to X]\mapsto [f_!\Q^H_Z]\:,$$
which commutes with pushdown $f_!$, exterior product $\boxtimes$ and pullback $g^*$. 
For $X=pt$ this corresponds to the ring homomorphism (\ref{H-ring}) under the identification
of $MHM(pt)\simeq mHs^p$.
\end{cor}

Here $K_0(var/X)$ is the motivic {\em relative Grothendieck group} of complex algebraic varieties over $X$,
i.e. the free abelian group generated by isomorphism classes $[f]=[f: Z\to X]$ of morphisms $f$ to $X$,
divided out be the {\em additivity relation} 
$$[f]=[f\circ i] + [f\circ j]$$ 
for a closed inclusion
$i: Y\to Z$ with open complement $j: U=Z\backslash Y\to Z$. The  pushdown $f_!$, exterior product $\boxtimes$ and pullback $g^*$
for these relative Grothendieck groups 
are defined by composition, exterior product and pullback of arrows. The fact that $\chi_{Hdg}$ commutes
with exterior product $\boxtimes$ (or pullback $g^*$) follows then from the corresponding K\"{u}nneth (or base change) theorem
for the functor 
$$f_!: D^bMHM(Z)\to D^bMHM(X)$$ 
(contained in Saito's work \cite{Sa4} and \cite{Sa1}[(4.4.3)]).\\

Let $\Lef:=[\A^1_{\C}]\in K_0(var/pt)$ be the class of the affine line so that
$$\chi_{Hdg}(\Lef)=[H^2(P^1(\C),\Q)]=[ \Q(-1)]\in K_0(MHM(pt))=K_0(mHs^p)$$
is the Lefschetz class $[ \Q(-1)]$. This is invertible in $K_0(MHM(pt))=K_0(mHs^p)$
so that the transformation $\chi_{Hdg}$ of corollary \ref{chiHdg} factorizes over the localization 
$$M_0(var/X):=K_0(var/X)[\Lef^{-1}] \:.$$
Altogether we get the following diagram of natural transformations commuting with
$f_!$, $\boxtimes$ and $g^*$:
\begin{equation}\label{motfunct} \begin{CD}
F(X) @< can << M_0(var/X) @<<< K_0(var/X) \\
@A \chi_{stalk} AA @VV \chi_{Hdg}V \\
K_0(D^b_c(X)) @<< rat < K_0(MHM(X)) \:. 
\end{CD}\end{equation}

Here $F(X)$ is the group of algebraically constructible functions on $X$ generated by
$1_Z$ for $Z\subset X$ a closed complex algebraic subset, with $ \chi_{stalk}$ given by the
Euler characteristic of the stalk complexes (compare \cite{Sc}[sec.2.3]). The pushdown
$f_!$ for algebraically constructible functions is defined for a morphism $f: Y\to X$ by
$$f_!(1_Z)(x):=\chi\left(H^*_c(Z\cap \{f=x\},\Q)\right) \quad \text{for $x\in X$,}$$
so that the horizontal arrow $can$ is given by
$$can: [f: Y\to X]\mapsto f_!(1_Y) \:, \quad \text{with $can(\Lef)=1_{pt}$.}$$ 

The advantage of  $M_0(var/X)$ compared to $K_0(var/X)$ is the fact, that it has an induced 
{\em duality} involution  $\DC: M_0(var/X)\to M_0(var/X)$ characterized uniquely by
(compare \cite{Bi}):
$$\DC\left([f: M\to X]\right)=\Lef^{-m}\cdot [f: M\to X]$$
for $f: M\to X$ a proper morphism with $M$ smooth and pure $m$-dimensional. This ``motivic duality" $\DC$ commutes with pushdown $f_!$ for proper $f$, so that $\chi_{Hdg}$ also commutes with duality by
\begin{equation}\label{motduality}\begin{split}
\chi_{Hdg}\left(\DC[id_M]\right) &= \chi_{Hdg}\left(\Lef^{-m}\cdot [Id_M]\right)
=[\Q^H_M(m)]\\
&= [\Q^H_M[2m](m)]= [\DC( \Q^H_M)]= \DC\left(\chi_{Hdg}\left([Id_M]\right)\right)
\end{split}\end{equation}
for $M$ smooth and pure $m$-dimensional. In fact by resolution of singularities and ``additivity", $K_0(var/X)$ is generated by such classes $f_![id_M]=[f: M\to X]$.\\

Then all the transformations of (\ref{motfunct}) {\em commute with duality}, were 
$K_0(D^b_c(X))$ gets this involution from  Verdier duality, and $\DC=id$ for
algebraically constructible functions by $\chi\left([ \Q(-1)]\right)=1_{pt}$ (compare also with \cite{Sc}[sec.6.0.6]). Similarly they commute with $f_*$ and $g^!$ defined by the relations (compare \cite{Bi}):
$$\DC\circ g^* = g^!\circ \DC \quad \text{and} \quad
\DC\circ f_* = f_!\circ \DC \:.$$
For example for an open inclusion $j: M\to \bar{M}$ one gets
\begin{equation}\label{j-dual}
\chi_{Hdg}\left(j_*[id_M]\right)=j_*[\Q^H_M] \:.
\end{equation}

\section{Characteristic classes of mixed Hodge modules}
\subsection{Homological characteristic classes}
In this section we explain the theory of $K$-theoretical characteristic homology classes of mixed Hodge modules based on the following result of Saito
(compare with \cite{Sa0}[sec.2.3] and \cite{Sa5}[sec.1] for the first part, and with \cite{Sa1}[sec.3.10, prop.3.11]) for the part (2)):

\begin{thm}[M. Saito]\label{grDR}
Let $Z$ be a complex algebraic variety. Then there is a functor of triangulated categories
\begin{equation} Gr^F_pDR: D^bMHM(Z) \to D^b_{coh}(Z)\end{equation}
commuting with proper push-down,
with $Gr^F_pDR(\M)=0$ for almost all $p$ and $\M$ fixed,
where $D^b_{coh}(Z)$ is the bounded
derived category of sheaves of algebraic $\OO_Z$-modules with coherent
cohomology sheaves. If $M$ is a (pure $m$-dimensional) complex algebraic manifold, then one
has in addition:
\begin{enumerate}
\item Let $\M\in MHM(M)$ be a single mixed Hodge module. Then
$Gr^F_pDR(\M)$ is the corresponding complex associated to the de Rham complex of the underlying algebraic left $D$-module $\M$ with its integrable connection $\nabla$:
$$ \begin{CD} DR(\M)=[\M @> \nabla >> \cdots @> \nabla >> \M\otimes_{\OO_M} \Omega^m_M]
\end{CD}$$ 
with $\M$ in degree $-m$, filtered  by
$$ \begin{CD}
F_p DR(\M) =[F_p\M @> \nabla >> \cdots  @> \nabla >> F_{p+m}\M\otimes_{\OO_M} \Omega^m_M] \:.
\end{CD} $$
\item Let $\bar{M}$ be a smooth partial compactification of the complex algebraic manifold $M$ with complement $D$ a normal crossing divisor with smooth irreducible components, with $j: M\to \bar{M}$ the open inclusion. 
Let $\VB=(L,F,W)$ be a good" variation of mixed Hodge structures on $M$. 
Then the filtered de Rham complex 
$$(DR(j_*\VB),F) \quad \text{of}\quad  j_*\VB\in MHM(\bar{M})[-m]\subset D^bMHM(\bar{M})$$
is filtered quasi-isomorphic to the logarithmic de Rham complex $DR_{log}(\LL)$ with the increasing filtration $F_{-p}:=F^{p}$ ($p\in \Z$) associated to the decreasing $F$-filtration
(\ref{drlog-fil}). In particular $Gr^F_{-p}DR(j_*\VB)$ ($p\in \Z$) is quasi-isomorphic to
$$\begin{CD}
Gr_F^{p}DR_{log}\left(\overline{\LL}\right) =
[Gr_F^{p}\overline{\LL} @> Gr\;\overline{\nabla} >> \cdots @> Gr\;\overline{\nabla} >> Gr_F^{p-m}\overline{\LL}\otimes_{\OO_{\bar{M}}} \Omega^m_{\bar{M}}(log(D))] \:.
\end{CD}$$
\end{enumerate}
\end{thm}

Here the filtration $F_p DR(\M)$ of the de Rham complex is well defined, since 
the action of the integrable connection $\nabla$ is given in local coordinates $(z_1,\dots, z_m)$ by
$$\nabla(\cdot) = \sum_{i=1}^m\: \frac{\partial}{\partial z_i}(\cdot) \otimes dz_i \:,
\quad \text{with
$ \frac{\partial}{\partial z_i}\in F_1\DC_M$,}$$
so that $\nabla(F_p\M)\subset
F_{p+1}\M$ for all $p$ by (\ref{D-filt}).
For later use, let us point that the maps $Gr\;\nabla$ and
$Gr\;\overline{\nabla}$ in the complexes 
$$Gr^F_{p}DR(\M) \quad \text{and} \quad
Gr_F^{p}DR_{log}\left(\overline{\LL}\right)$$ 
are $\OO$-linear!

\begin{example}\label{vanishing}
Let $M$ be a pure $m$-dimensional complex algebraic manifold. Then
$$Gr^F_{-p}DR(\Q^H_M)\simeq \Omega^p_M[-p] \in D^b_{coh}(M)$$
for $0\leq p \leq m$, and $Gr^F_{-p}DR(\Q^H_M)\simeq 0$ otherwise.
Assume in addition that $f:M\to Y$ is a resolution of singularities of the
pure dimensional complex algebraic variety $Y$. Then 
$IC^H_Y$ is a direct summand of $f_*\Q^H_M \in D^bMHM(Y)$ so that by functoriality
$gr^F_{-p}DR(IC^H_Y)$ is a direct summand of $Rf_*\Omega^p_M[-p] \in D^b_{coh}(Y)$.
In particular $$Gr^F_{-p}DR(IC^H_Y)\simeq 0 \quad \text{ for $p<0$ or $p>m$.}$$
\end{example}

The transformations
$Gr^F_pDR$ ($p\in \Z$) induce functors on the level of Grothendieck groups.
Therefore, if $G_0(Z) \simeq K_0(D^b_{coh}(Z))$ denotes the
Grothendieck group of coherent {\em algebraic} $\OO_Z$-sheaves on $Z$, we get group
homomorphisms 
$$Gr^F_pDR: K_0(MHM(Z))=K_0(D^bMHM(Z))\to K_0(D^b_{coh}(Z))\simeq G_0(Z)\:.$$

\begin{defn}
The {\em motivic Hodge Chern class transformation}
$$MHC_y: K_0(MHM(Z)) \to G_0(Z) \otimes \Z[y^{\pm 1}]$$
is defined by
\begin{equation}\label{MHC}
[\M] \mapsto \sum_{i,p}\; (-1)^{i} [\HC^i ( Gr^F_{-p} DR(\M) )] \cdot
(-y)^p\:.
\end{equation}
\end{defn}

So this characteristic class captures information from the graded pieces of the filtered
de Rham complex of the filtered $D$-module underlying a mixed Hodge module $\M\in MHM(Z)$, instead of the graded pieces of the filtered $D$-module itself (as more often studied). Let 
$p'=min\{p|\;F_p\M \neq 0\}$. Using  theorem \ref{grDR}(1) for a local embedding
$Z\hookrightarrow M$ of $Z$ into a complex algebraic manifold $M$ of dimension $m$,
one gets 
$$Gr^F_{p} DR(\M)=0 \quad \text{for $p<p'-m$, and} \quad 
Gr^F_{p'-m} DR(\M) \simeq \left(F_{p'}\M\right) \otimes_{\OO_M} \omega_M $$
is a coherent $\OO_Z$-sheaf independent of the local embedding. Here we are using left
 $D$-modules (related to variation of Hodge structures), whereas for this question the corresponding filtered right $D$-module (as used in \cite{Sa3})
$$\M^r:= \M\otimes_{\OO_M} \omega_M  \quad \text{with} \quad
F_p\M^r:=\left(F_{p+m}\M\right) \otimes_{\OO_M} \omega_M$$ 
would better work. Then the coefficient of the ``top-dimensional" power of $y$ in 
$MHC_y\left([\M]\right)$:
\begin{equation} 
MHC_y\left([\M]\right)= [F_{p'}\M \otimes_{\OO_M} \omega_M ]\otimes (-y)^{m-p'}
+ \sum_{i<m-p'} (\cdots) \cdot y^i \in G_0(Z)[y^{\pm 1}]
\end{equation}
is given by the class $ [F_{p'}\M \otimes_{\OO_M} \omega_M ] \in G_0(Z)$ of this 
coherent $\OO_Z$-sheaf (up to a sign). Using resolution of singularities, one gets for example
for an $m$-dimensional complex algebraic variety $Z$, that 
$$MHC_y([\Q^H_Z]) = [\pi_*\omega_M] \cdot y^m + \sum_{i<m} (\cdots) \cdot y^i \in G_0(Z)[y^{\pm 1}]\:,$$
with $\pi: M\to Z$ any resolution of singularities of $Z$ (compare \cite{Sa5}[cor.0.3]).
More generally, for an irreducible complex variety $Z$ and $\M=IC^H_Z(\LL)$ a pure Hodge module
with strict support $Z$, the corresponding coherent $\OO_Z$-sheaf
$$S_Z(\LL):=F_{p'}IC^H_Z(\LL) \otimes_{\OO_M} \omega_M$$
only depends on $Z$ and the good variation of Hodge structures $\LL$ on a Zariski open smooth subset of $Z$, and it behaves  much like a  dualizing sheaf. 
Its formal properties are studied in 
Saito's proof  given in  \cite{Sa3} of a conjecture of Kollar. So the ``top-dimensional" power of $y$ in  $MHC_y\left([IC^H_Z(\LL)]\right)$ exactly picks out (up to a sign) the class
$[S_Z(\LL)]\in G_0(Z)$ of this interesting coherent sheaf $S_Z(\LL)$ on $Z$.\\

Let $td_{(1+y)}$ be the {\em twisted Todd transformation}
\begin{equation}\label{td} \begin{split}
td_{(1+y)}:\: &G_0(Z) \otimes \Z[y^{\pm 1}] \to H_*(Z) \otimes
\Q[y^{\pm 1}, (1+y)^{-1}]\;;\\
&[\FC] \mapsto \sum_{k \geq 0}\; td_k([\FC]) \cdot (1+y)^{-k}\:,
\end{split} \end{equation}
where $H_*(\cdot)$
stands either for Chow homology groups $CH_*(\cdot)$ or for Borel-Moore homology groups $H_{2*}^{BM}(\cdot)$ (in even degrees), and $td_k$ is the
degree $k$ component in $H_{k}(Z)$ of the {\em Todd class
transformation} $td_*: G_0(Z) \to H_{*}(Z) \otimes \Q$ of
Baum-Fulton-MacPherson \cite{BFM}, which is linearly extended over
$\Z[y^{\pm 1}]$ (compare also with \cite{Fu}[chap.18] and \cite{FM}[Part II]).

\begin{defn}
The (un)normalized {\em motivic Hirzebruch class
transformations} $MHT_{y*}$ (and $MH\tilde{T}_{y*}$) are defined by the composition
\begin{equation}\label{MHT}
MHT_{y*} :=td_{(1+y)} \circ MHC_y: K_0(MHM(Z)) \to H_{*}(Z)
\otimes \Q[y^{\pm 1},(1+y)^{-1}]
\end{equation}
and
\begin{equation}\label{MHT2}
MH{T}_{y*} :=td_* \circ MHC_y: K_0(MHM(Z)) \to H_{*}(Z)
\otimes \Q[y^{\pm 1}] \:.
\end{equation}
\end{defn}

\begin{remark}
By precomposing with the transformation $\chi_{Hdg}$ from corollary \ref{chiHdg}
one gets similar transformations
$$mC_y:=MHC_y\circ \chi_{Hdg},\; T_{y*}:=MHT_{y*}\circ \chi_{Hdg} \quad
\text{and} \quad \tilde{T}_{y*}:=MH\tilde{T}_{y*}\circ \chi_{Hdg}$$
defined on the relative Grothendieck group of complex algebraic varieties $K_0(Var/\cdot)$
as studied in
\cite{BSY}. Then it is the (normalized) motivic Hirzebruch class transformation $T_{y*}$, which ``unifies" in a functorial way 
\begin{enumerate}
\item[(-1)] the (rationalized) Chern class transformation $c_*$ of MacPherson \cite{M}, 
\item[(0)] the Todd class transformation $td_*$ of
Baum-Fulton-MacPherson \cite{BFM}, and 
\item[(1)] the $L$-class  transformation $L_*$ 
of Cappell-Shaneson \cite{CS} 
\end{enumerate}
for $y=-1,0$ and $1$ respectively (compare with \cite{BSY, SY} and also with \cite{Y} in these proceedings).
\end{remark}

In this paper we work most the time only with the more important $K$-theoretical
transformation $MHC_y$. The corresponding results for $MHT_{y*}$ follow from this by the known
properties of the Todd class transformation $td_*$ (compare \cite{BFM, Fu, FM}).

\begin{example}\label{pt}\rm Let $\VB=(V,F,W) \in MHM(pt)=mHs^p$ be a (graded polarizable) mixed Hodge structure. Then:
\begin{equation}\label{point}
MHC_y([\VB])= \sum_p\:
\text{dim}_{\C} (Gr^p_F V_{\C}) \cdot (-y)^p = \chi_y([\VB]) \in \Z[y^{\pm 1}]
=G_0(pt)\otimes \Z[y^{\pm 1}]\:.
\end{equation}
So over a point the transformation $MHC_y$ coincides with the
$\chi_y$-genus ring homomorphism $\chi_y:K_0(mHs^p) \to
\Z[y^{\pm 1}]$ (and similarly for $MH\tilde{T}_{y*}$ and $MHT_{y*}$).
\end{example}

The {\em motivic Chern resp. Hirzebruch class} $C_y(Z)$ resp. $T_{y*}(Z)$ of a complex algebraic
variety $Z$ is defined by \begin{equation}
C_y(Z):=MHC_y([\Q_Z^H]) \quad \text{and} \quad
T_{y*}(Z):=MHT_{y*}([\Q_Z^H]) \:.
\end{equation}
Similarly, if $U$ is a pure $n$-dimensional complex algebraic manifold,
and $\LL$ is a local system on $U$ underlying a good
variation of mixed Hodge structures, we define {\em twisted motivic Chern resp. Hirzebruch
characteristic classes} by (compare \cite{CLMS, CLMS2, MSc})
\begin{equation}\label{tHc} 
C_y(U; \LL):=MHC_y([\LL^H]) \quad \text{and} \quad
T_{y*}(U; \LL):=MHT_{y*}([\LL^H]) \:,\end{equation} 
where $\LL^H[n]$ is the smooth mixed Hodge module on $U$ with underlying
perverse sheaf $\LL[n]$. Assume in addition, that $U$ is dense and Zariski open in the complex algebraic variety $Z$. Let $IC^H_Z, IC_Z^H(\LL) \in MHM(Z)$ be the (twisted) intersection
homology (mixed) Hodge module on $Z$, whose underlying perverse sheaf is $IC_Z$ resp.
$IC_Z(\LL)$.
Then we define {\em Intersection characteristic classes} by (compare \cite{BSY, CMS, CLMS2, MSc}):
\begin{equation} IC_y(Z):=MHC_y\left(\left[IC^H_Z[-n]\right]\right) \quad \text{and} \quad
IT_{y*}(Z):=MHT_{y*}\left(\left[ IC^H_Z[-n]\right]\right) \end{equation} 
and similarly,
\begin{equation}
IC_y(Z; \LL):=MHC_y\left(\left[IC_Z^H(\LL)[-n] \right]\right) \:\: \text{and} \:\:
IT_{y*}(Z;\LL):=MHT_{y*}\left(\left[ IC^H_Z(\LL)[-n]\right]\right) \:.\end{equation}

By definition and theorem \ref{grDR}, the transformations $MHC_y$ and $MHT_{y*}$   
{\em commute with proper
push-forward}. The following {\it normalization} property holds (cf.
\cite{BSY}): If $M$ is smooth, then
\begin{equation}
C_y(Z)=\lambda_y(T^*M)\cap [\OO_M] \quad \text{and} \quad
T_{y*}(Z)=T_y^*(TM) \cap [M]\:,
\end{equation}
where $T_y^*(TM)$ is the
cohomology Hirzebruch class of $M$ as in theorem \ref{gHRR}.

\begin{example}
Let $Z$ be a compact (possibly singular) complex algebraic  variety, with $k:Z \to pt$ the proper constant map to a point. Then for $\M\in D^bMHM(Z)$ the pushdown 
$$k_*(MHC_y(\M))=MHC_y(k_*\M)= \chi_y\left([H^*(Z,\M)]\right)$$ 
is the Hodge genus 
\begin{equation}\label{defchi} \chi_y([H^*(Z,\M)])=\sum_{i,p}\;
(-1)^i dim_{\C} (Gr^p_F H^i(Z,\M)) \cdot (-y)^p \:.
\end{equation}
In particular:
\begin{enumerate}
\item If $Z$ is smooth, then
$$k_*C_y(Z)=\chi_y(Z):=\chi_y\left([H^*(Z,\Q)]\right)$$
and
$$k_*C_y(Z;\LL)=\chi_y(Z;\LL):=\chi_y\left([H^*(Z,\LL)]\right)\:.$$
\item If $Z$ is pure dimensional, then
$$k_*IC_y(Z)=I\chi_y(Z):=\chi_y\left([IH^*(Z,\Q)]\right)$$
and
$$k_*IC_y(Z;\LL)=I\chi_y(Z;\LL):=\chi_y\left([IH^*(Z,\LL)]\right)\:.$$
\end{enumerate}
\end{example}

Note that for $Z$ compact
$$I\chi_{-1}(Z)=\chi([IH^*(Z;\Q)]$$  is the
{\em intersection (co)homology Euler characteristic} of $Z$, whereas for $Z$
projective, 
$$I\chi_1(Z)=sign\left(IH^*(Z,\Q)\right)$$ 
is the {\em intersection (co)homology signature} of $Z$ due to Goresky-MacPherson \cite{GM2}. In fact this follows as in the smooth context from
Saito's (relative version of the) Hodge index theorem for intersection cohomology (\cite{Sa0}[thm.5.3.2]).
Finally $\chi_0(Z)$ and $I\chi_0(Z)$ are two possible extensions to singular varieties of the
{\em arithmetic genus}. Here it makes sense to take $y=0$, since
one has by Example \ref{vanishing}:
$$k_*IC_y(Z)=I\chi_y(Z) \in \Z[y]\:.$$

It is conjectured that for a
pure $n$-dimensional  compact variety $Z$:
$${IT_1}_*(Z)\stackrel{?}{=} L_*(Z) \in H_{2*}(Z,\Q)$$
is the Goresky-MacPherson homology
$L$-class \cite{GM2} of the Witt space $Z$ (\cite{BSY}, Remark 5.4). Similarly one should expect for  a pure-dimensional compact variety $Z$, that 
\begin{equation}
\alpha(IC_1(Z)) \stackrel{?}{=} \triangle(Z) \in KO^{top}_0(Z)[1/2]\oplus  KO^{top}_2(Z)[1/2]\simeq K^{top}_0(Z)[1/2]\:,
\end{equation}
where $\alpha: G_0(Z)\to K_0^{top}(Z)$ is the K-theoretical Riemann-Roch transformation
of Baum-Fulton-MacPherson \cite{BFM2}, and
$\triangle(Z)$ is the {\em Sullivan class} of the Witt space $Z$ (compare with \cite{Ba2}
in these proceedings).
These conjectured equalities are true for a smooth $Z$, or more generally for a pure $n$-dimensional compact
complex algebraic variety $Z$ with a {\em small resolution} of singularities
$f: M\to Z$, in which case one has $f_*(\Q_M^H)=IC^H_Z[-n]$ so that
$$IT_{1*}(Z) = f_*T_{1*}(M) = f_*L_*(M)= L_*(Z) \:,$$
and
$$\alpha\left(IC_1(Z)\right) = f_*\left(\alpha(C_1(M))\right) = f_*\triangle(M)= \triangle(Z)\:.$$
Here the functoriality $f_*L_*(M)= L_*(Z)$ and $f_*\triangle(M)= \triangle(Z)$ for a small resolution
follows e.g. from the work \cite{Woo}, which allows one to think of the characteristic classes $L_*$ and $\triangle$
as covariant functors for suitable Witt groups of selfdual constructible sheaf complexes. \\

In particular, the classes $f_*C_1(M)$ and $f_*T_{1*}(M)$ do not depend on the choice of a small
resolution. In fact the same functoriality argument applies to (compare \cite{CMS, MSc})
$$IC_y(Z)=f_*C_y(M) \in G_0(Z)\otimes\Z[y] \quad \text{and} \quad IT_{y*}(Z)=f_*T_{y*}(M)\in 
H_{2*}(Z)\otimes \Q[y,(1+y)^{-1}] \:.$$
Note that in general a complex variety $Z$ doesn't have a small resolution,
and even if it exists, it is in general not unique.
This type of independence question were discussed by Totaro \cite{To},
pointing out the relation to the famous {\em elliptic genus and classes}
(compare also with \cite{Lib,Wae} in these proceedings).
Note that we  get such a result for the K-theoretical class 
$$IC_y(Z)=f_*C_y(M) \in G_0(Z)\otimes\Z[y] \:!$$

\subsection{Calculus of characteristic classes}
So far we only discussed the functoriality of $MHC_y$ with respect to proper push down,
and the corresponding relation to Hodge genera for compact $Z$ coming from the push down
for the proper constant map $k: Z\to pt$. Now we explain some other important
functoriality properties. Their proof is based on the following (e.g. see \cite{MSc}[(4.6)]):

\begin{example}\label{MHMopen}
Let $\bar{M}$ be a smooth partial compactification of the complex algebraic manifold $M$ with complement $D$ a normal crossing divisor with smooth irreducible components, with $j: M\to \bar{M}$ the open inclusion. 
Let $\VB=(L,F,W)$ be a good variation of mixed Hodge structures on $M$. 
Then the filtered de Rham complex 
$$(DR(j_*\VB),F) \quad \text{of}\quad  j_*\VB\in MHM(\bar{M})[-m]\subset D^bMHM(\bar{M})$$
is by theorem \ref{grDR}(2) filtered quasi-isomorphic to the logarithmic de Rham complex $DR_{log}(\LL)$ with the increasing filtration $F_{-p}:=F^{p}$ ($p\in \Z$) associated to the decreasing $F$-filtration (\ref{drlog-fil}). Then
\begin{equation}\begin{split}\label{MHCj*-fomula}
MHC_y(j_*\VB) &= \sum_{i,p}\; (-1)^{i} [\HC^i ( Gr_F^{p} DR_{log}(\LL) )] \cdot (-y)^p\\
&=\sum_{p}\; [Gr_F^{p} DR_{log}(\LL) ] \cdot (-y)^p\\
&\stackrel{(\ast)}{=} \sum_{i,p}\; (-1)^{i} [ Gr_F^{p-i}(\LL) \otimes_{\OO_{\bar{M}}} \Omega_{\bar{M}}^i(log(D))] \cdot (-y)^p\\
&=MHC^y(Rj_*L) \cap \left(\lambda_y\left( \Omega_{\bar{M}}^1(log(D))\right)\cap [\OO_{\bar{M}}] \right)\:.
\end{split}\end{equation}
In particular for $j=id: M\to M$ we get the following {\em Atiyah-Meyer type formula}
(compare \cite{CLMS, CLMS2, MSc}):
\begin{equation}\label{MHCy-fomula}
MHC_y(\VB) = MHC^y(L) \cap \left(\lambda_y(T^*M)\cap [\OO_M] \right)\:.
\end{equation}
\end{example}

\begin{remark}
The formula (\ref{MHCj*-fomula}) is a class version of the formula (\ref{pairing1}) of theorem \ref{j*VMHS}, which one gets back from (\ref{MHCj*-fomula}) by pushing down  to a point
for the proper constant map $k: \bar{M}\to pt$ on the compactification $\bar{M}$ of $M$.

Also note that in the equality ($\ast$) above we use the fact that the complex $Gr_F^{p} DR_{log}(\LL)$ has coherent (locally free) objects, with $\OO_{\bar{M}}$-linear maps between them.
\end{remark}

The formula (\ref{MHCj*-fomula}) describes a {\em splitting} of the characteristic class
$MHC_y(j_*\VB)$ into 
\begin{enumerate}
\item[(coh)] a cohomological term $MHC^y(Rj_*L)$ capturing the information of the
good variation of mixed Hodge structures $L$, and 
\item[(hom)] the homological term
$$\lambda_y\left( \Omega_{\bar{M}}^1(log(D))\right)\cap [\OO_{\bar{M}}] = MHC_y(j_*\Q^H_M)$$
capturing the information of the underlying space or embedding $j: M\to \bar{M}$.
\end{enumerate}

The term $MHC^y(Rj_*L)$ has by corollary \ref{Cor-j*} good functorial behavior with respect to exterior and suitable tensor products, as well as for smooth pullbacks. 
For the exterior products one gets similarly (compare \cite{De}[prop.3.2]):
$$\Omega_{\bar{M}\times \bar{M}'}^1(log(D\times M'\cup M\times D')) \simeq \left(\Omega_{\bar{M}}^1(log(D))\right) \boxtimes
\left( \Omega_{\bar{M'}}^1(log(D'))\right)$$
so that
$$\lambda_y\left( \Omega_{\bar{M}\times \bar{M}'}^1(log(D\times M'\cup M\times D'))\right)\cap 
[\OO_{\bar{M}\times \bar{M}'}] =$$
$$\left(\lambda_y\left( \Omega_{\bar{M}}^1(log(D))\right)\cap [\OO_{\bar{M}}] \right)\boxtimes
\left(\lambda_y\left( \Omega_{\bar{M'}}^1(log(D'))\right)\cap [\OO_{\bar{M'}}] \right)
$$
for the product of two partial compactifications as in example \ref{MHMopen}.
But the Grothendieck group $K_0(MHM(Z))$ of mixed Hodge modules on the complex variety $Z$ is generated by classes of the form $f_*(j_*[\VB])$, with $f: \bar{M}\to Z$ proper and $M,\bar{M},\VB$ as
before. Finally one also has the multiplicativity
$$(f\times f')_*= f_*\boxtimes f'_*$$
for the push down for proper maps $f: \bar{M}\to Z$ and  $f': \bar{M}'\to Z'$
on the level of Grothendieck groups $K_0(MHM(\cdot))$ as well as for $G_0(\cdot)\otimes\Z[y^{\pm 1}]$. Then one gets (as in \cite{BSY}[Proof of Cor. 2.1(3)]) from corollary \ref{Cor-j*}
and the example \ref{MHMopen} the following

\begin{cor}[Multiplicativity for exterior products] The motivic Chern class transformation
$MHC_y$ commutes with exterior products:
\begin{equation}
MHC_y([M\boxtimes M'])=
MHC_y([M]\boxtimes [M'])= MHC_y([M]) \boxtimes MHC_y([M'])
\end{equation}
for $M\in D^bMHM(Z)$ and $M'\in D^bMHM(Z')$.
\end{cor}

Next we explain the behaviour of $MHC_y$ for smooth pullbacks. Consider a cartesian
diagram of morphisms of complex algebraic varieties
$$\begin{CD}
\bar{M}' @> g' >> \bar{M}\\
@V f' VV @VV f V\\
Z' @>> g > Z \:,
\end{CD}$$
with $g$ smooth, $f$ proper and $M,\bar{M},\VB$ as before.
Then also $g'$ is  smooth and $f'$ is proper, and one has the {\em base change isomorphism}
$$g^*f_*=f'_*g'^*$$
on the level of Grothendieck groups $K_0(MHM(\cdot))$ as well as for $G_0(\cdot)\otimes\Z[y^{\pm 1}]$.
Finally for the induced partial compactification $\bar{M}'$ of $M':=g'^{-1}(M)$, with complement
$D'$ the induced normal crossing divisor with smooth irreducible components,
one has a short exact sequence of vector bundles on $\bar{M}'$:
$$0\to g'^*\left(\Omega_{\bar{M}}^1(log(D))\right)\to \Omega_{\bar{M}'}^1(log(D'))
\to T_{g'}^* \to 0\:,$$
with $T_{g'}^*$ the relative cotangent bundle along the fibers of the smooth morphism $g'$.
And by base change one has $T_{g'}^*=f'^*(T_g^*)$. So for the corresponding lambda classes we get
\begin{equation}\begin{split}
\lambda_y\left( \Omega_{\bar{M}'}^1(log(D'))\right) &=
\left(g'^*\lambda_y\left( \Omega_{\bar{M}}^1(log(D))\right)\right)\otimes \lambda_y(T_{g'}^*)\\
&=\left(g'^*\lambda_y\left( \Omega_{\bar{M}}^1(log(D))\right)\right)\otimes
f'^*\lambda_y(T_g^*) \:.
\end{split}\end{equation}

Finally (compare also with \cite{BSY}[Proof of Cor. 2.1(4)]), by using the {\em projection formula}
$$\lambda_y(T_g^*) \otimes f'_*(\cdot)= f'_*\left( f'^*\lambda_y(T_g^*)\otimes (\cdot)\right)\; : \:,\:\: G_0(\bar{M}')\otimes\Z[y^{\pm 1}]\to G_0(Z')\otimes\Z[y^{\pm 1}]
$$
one gets from corollary \ref{Cor-j*}
and the example \ref{MHMopen} the following

\begin{cor}[VRR for smooth pullbacks] \label{VRR}
For a smooth morphism $g: Z'\to Z$ of complex algebraic varieties one has for
the motivic Chern class transformation the following {\em Verdier Riemann-Roch formula}:
\begin{equation}
\lambda_y(T_g^*) \cap 
g^*MHC_y([M])= MHC_y(g^*[M])=  MHC_y([g^*M])
\end{equation}
for $M\in D^bMHM(Z)$. In particular
\begin{equation}
g^*MHC_y([M])= MHC_y(g^*[M])=  MHC_y([g^*M])
\end{equation}
for $g$ an \'{e}tale morphism (i.e. a smooth morphism with zero dimensional fibers), or in more topological terms, for $g$ an unramified covering.
The most important special case is that of an open embedding.
\end{cor}

If moreover $g$ is also proper, then one gets from corollary \ref{VRR} and the projection formula the following result:

\begin{cor}[Going up and down] \label{updown}
Let $g: Z'\to Z$ be a smooth and proper  morphism of complex algebraic varieties.
Then one has for
the motivic Chern class transformation the following {\em going up und down formula}:
\begin{equation}\begin{split}
MHC_y(g_*g^*[M])&=g_*MHC_y(g^*[M])\\
&=g_*\left( \lambda_y(T_g^*) \cap g^*MHC_y([M])\right)\\
&=\left(g_*\lambda_y(T_g^*)\right) \cap  MHC_y([M])  
\end{split}\end{equation}
for $M\in D^bMHM(Z)$, with 
$$g_*\left(\lambda_y(T_g^*)\right):= \sum_{p,q\geq 0}\;(-1)^q\cdot [R^qg_*(\Omega^p_{Z'/Z})]\cdot y^p
\in K^0_{alg}(Z)[y]$$
the algebraic cohomology class being given (as in example \ref{ex:smooth}) by 
$$MHC^y([Rg_*\Q_{Z'}])= \sum_{p,q\geq 0}\;(-1)^q\cdot [R^qg_*(\Omega^p_{Z'/Z})]\cdot y^p\:.$$
Note that all higher direct image sheaves $R^qg_*(\Omega^p_{Z'/Z})$ are locally free in this case, since $g$ is a smooth and proper morphism of complex algebraic varieties
(compare with \cite{De2}). In particular
$$g_*C_y(Z')= \left(g_*\lambda_y(T_g^*)\right) \cap C_y(Z) \:,$$
and 
$$g_*IC_y(Z')= \left(g_*\lambda_y(T_g^*)\right) \cap IC_y(Z) $$
for $Z$ and $Z'$ pure dimensional. If moreover $Z,Z'$ are compact, with $k: Z\to pt$ the constant proper map, then
\begin{equation}
\chi_y(g^*[\M])=k_*g_*MHC_y(g^*[M])=\langle g_*\lambda_y(T_g^*),  MHC_y([M])\rangle \:. 
\end{equation}
In particular
$$\chi_y(Z')= \langle g_*\lambda_y(T_g^*), C_y(Z)\rangle \quad \text{and} \quad
I\chi_y(Z')= \langle g_*\lambda_y(T_g^*), IC_y(Z)\rangle \:.$$
\end{cor}

The result of this corollary can also be seen form a different view point, by making the
``going up and down" calculation already on the level of Grothendieck groups of mixed Hodge modules, where this time one only needs the assumption that $f: Z'\to Z$ is proper
(to get the projection formula):
$$f_*f^*[\M]=[f_*f^*\M]=[f_*(\Q^H_{Z'}\otimes f^*\M)] = [f_* \Q^H_{Z'}]\otimes [\M] \in K_0(MHM(Z))$$
for $\M\in D^bMHM(Z)$. The problem for a singular $Z$ is then, that we do not have a precise
relation between 
$$[f_* \Q^H_{Z'}] \in K_0(MHM(Z)) \quad \text{and} \quad [Rf_* \Q_{Z'}]\in K_0(FmHs^p(Z))\:.$$

\begin{remark}
What is {\em missing} up to now is the right notion of a good variation (or family) of mixed
Hodge structures on a {\em singular} complex algebraic variety $Z$! 
This class should at least contain:
\begin{enumerate}
\item The higher direct image local systems $R^i f_* \Q_{Z'}$ ($i\in \Z$) for a smooth and proper morphism
$f: Z'\to Z$ of complex algebraic varieties.
\item The pullback $g^*\LL$ of a good variation of mixed Hodge structures $\LL$ on a smooth complex algebraic manifold $M$ under
an algebraic morphism $g: Z\to M$.
\end{enumerate}
At the moment we have to assume that $Z$ is smooth (and pure dimensional) so that one can use theorem \ref{MHM2}.
\end{remark}

Nevertheless, in case (2) above we can already prove the following interesting result
(compare with \cite{MSc}[sec.4.1] for a similar result for $MHT_{y*}$ in the case 
when $f$ is a closed embedding):

\begin{cor}[Multiplicativity]\label{pairing-class} Let $f: Z\to N$ be a morphism of complex algebraic varieties, with $N$ smooth and pure $n$-dimensional.
Then one has a natural {\em pairing}
$$f^*(\cdot) \cap (\cdot): K_0(VmHs^g(N)) \times K_0(MHM(Z)) \to K_0(MHM(Z))\:,$$
$$([\LL],[\M])\mapsto [f^*\left(\LL^H\right)\otimes \M] \:.$$
Here $\LL^H[m]$ is the smooth mixed Hodge module on $N$ with underlying perverse sheaf $\LL[m]$. 
One also has a similar pairing on (co)homological level:
$$f^*(\cdot) \cap (\cdot): K^0_{alg}(N)\otimes \Z[y^{\pm 1}] \times G_0(Z)\otimes \Z[y^{\pm 1}]  \to G_0(Z)\otimes \Z[y^{\pm 1}] \:,$$
$$([\VV]\cdot y^i,[\FC]\cdot y^j)\mapsto [f^*(\VV)\otimes \FC]\cdot y^{i+j} \:.$$
And the motivic Chern class transformations $MHC^y$ and
$MHC_y$ commute with these natural pairings:
\begin{equation}\begin{split}
MHC_y\left([f^*\left(\LL^H\right)\otimes \M]\right) &= MHC^y\left([f^*\LL]\right) \cap  MHC_y([\M])\\
&= f^*\left(MHC^y\left([\LL]\right)\right) \cap  MHC_y([\M])
\end{split}\end{equation}
for $\LL \in VmHs^g(N)$ and $\M\in D^bMHM(Z)$.
\end{cor}

For the proof 
we can once more assume $\M=g_*j_*\VB$ for $g: \bar{M}\to Z$ proper, with $\bar{M}$ a pure-dimensional smooth complex algebraic manifold,
$j: M\to \bar{M}$ a Zariski open inclusion with complement $D$ a normal crossing divisor with smooth irreducible components, and finally
$\VB$ a good variation of mixed Hodge structures on $M$. Using the projection formula, it is then enough to prove
$$MHC_y\left([g^*f^*\left(\LL^H\right)\otimes j_*\VB]\right) = MHC^y\left([g^*f^*\LL]\right) \cap  MHC_y([j_*\VB]) \:.$$
But $g^*f^*\LL$ is a good variation of mixed Hodge structures on $\bar{M}$, so that by example \ref{MHMopen} and corollary \ref{Cor-j*}(3)
both sides are equal to
$$\left(MHC^y(g^*f^*\LL) \otimes MHC^y(j_*\VB)\right) \cap  \left(\lambda_y\left( \Omega_{\bar{M}}^1(log(D))\right)\cap [\OO_{\bar{M}}] \right)\:.$$

As an application of the very special case $f=id: Z\to N$ the identity of a complex algebraic manifold $Z$, with
$$MHC_y([\Q^H_Z])= \lambda_y(T^*Z) \cap [\OO_Z] \:,$$
one gets the  Atiyah-Meyer type formula (\ref{MHCy-fomula}) as well as the following result
(cf. \cite{CLMS, CLMS2, MSc}):

\begin{example}[Atiyah type formula] Let $g: Z' \to Z$ be a proper morphism of complex algebraic varieties, with $Z$ smooth and connected. Assume that for a given
$\M\in D^bMHM(Z')$ all direct image sheaves
$$R^ig_*rat(\M)\quad \text{($i\in \Z)$ are locally constant}\:,$$ 
e.g. $g$ is a locally trivial fibration and $\M=\Q^H_{Z'}$ or $\M=IC^H_{Z'}$
(for $Z'$ pure dimensional), so that they all underlie a good variation of mixed
Hodge structures. Then one can define 
$$[Rg_*rat(\M)]:=\sum_{i\in \Z}\; (-1)^i\cdot [R^ig_*rat(\M)] \in K_0(VmHs^g(Z))\:,$$
with
\begin{equation}\begin{split}
g_*MHC_y([\M])&= MHC_y(g_*[M])\\
&= MHC^y([Rg_*rat(\M)])\otimes 
\left(\lambda_y(T^*Z) \cap [\OO_Z]\right) \:.
\end{split}\end{equation}
\end{example}

As a final application we mention the

\begin{example}[Atiyah-Meyer type formula for intersection cohomology] \label{AM-IC}
 Let $f: Z\to N$ be a morphism of complex algebraic varieties, with $N$ smooth and pure $n$-dimensional (e.g. a closed embedding).
Assume also $Z$ is pure $m$-dimensional. Then one has for a good variation of mixed Hodge structures 
$\LL$ on $N$ the equality
$$IC^H_Z(f^*\LL)[-m] \simeq f^*\LL^H \otimes IC^H_Z[-m] \in MHM(Z)[-m]\subset D^bMHM(Z)\:,$$
so that
\begin{equation}
IC_y(Z;f^*\LL) = MHC^y(f^*\LL) \cap IC_y(Z) = f^*\left(MHC^y(\LL)\right) \cap IC_y(Z) \:.
\end{equation}
If in addition $Z$ is also compact, then one gets by pushing down to a point:
\begin{equation}
 I\chi_y(Z;f^*\LL)= \langle MHC^y(f^*\LL), IC_y(Z)\rangle \:.
\end{equation}
\end{example}

\begin{remark}
This example should be seen as a Hodge theoretical version of the corresponding
result of Banagl-Cappell-Shaneson \cite{BCS} for the $L$-classes
$L_*(IC_Z(L))$ of a selfdual {\em Poincar\'{e} local system} $L$ on all of $Z$.
The special case of example \ref{AM-IC} for $f$ a closed inclusion was already explained
in \cite{MSc}[sec.4.1]. 
\end{remark}

Finally note that all the results of this section can easily be applied to the
(un)normalized {\em motivic Hirzebruch class transformation} $MHT_{y*}$ (and  $MH\tilde{T}_{y*}$), because the {\em Todd class transformation} $td_*: G_0(\cdot) \to
H_*(\cdot)$ of Baum-Fulton-MacPherson \cite{BFM} has the following properties
(compare also with \cite{Fu}[chapter 18] and \cite{FM}[Part II]):

\begin{description}
\item[Functoriality] The Todd class transformation $td_*$ commutes with pushdown 
$f_*$ for a proper morphism $f: Z\to X$: 
$$td_*\left(f_*\left([\FC]\right)\right)=
f_*\left(td_*\left([\FC]\right)\right) \quad \text{for $[\FC]\in G_0(Z)$.}$$
\item[Multiplicativity for exterior products] The  Todd class transformation $td_*$ commutes with exterior products: 
$$td_*\left([\FC\boxtimes \FC']\right)=
td_*\left([\FC]\right) \boxtimes td_*\left([\FC']\right) \quad \text{for
$[\FC]\in G_0(Z)$ and $[\FC']\in G_0(Z')$.}$$
\item[VRR for smooth pullbacks] For a smooth morphism $g: Z'\to Z$ of complex algebraic varieties one has for
the Todd class transformation $td_*$ the following  Verdier Riemann-Roch formula:
$$td^*(T_g) \cap 
g^*td_*([\FC])= td_*(g^*[\FC])= td_*([g^*\FC]) \quad \text{for $[\FC]\in G_0(Z)$.}$$
\item[Multiplicativity] Let $ch^*: K^0_{alg}(\cdot)\to H^*(\cdot)\otimes \Q$ be the cohomological {\em Chern character} to the cohomolgy $H^*(\cdot)$  given by the 
operational Chow ring $CH^*(\cdot)$ or the usual cohomology $H^{2*}(\cdot,\Z)$ in even degrees. Then one has the multiplicativity relation
$$td_*([\VV\otimes \FC])= ch^*([\VV]) \cap td_*([\FC]) $$
for $[\VV]\in  K^0_{alg}(Z)$ and $[\FC]\in G_0(Z)$, with $Z$ a (possible singular)
complex algebraic variety.
\end{description}

\subsection{Characteristic classes and duality}
In this final section we explain the characteristic class version of the duality formula
(\ref{chi-y-dual}) for the $\chi_y$-genus. We also show that the
specialization of $MHT_{y*}$ for $y=-1$ exists and is equal to the rationalized 
MacPherson Chern class $c_*$ of the underlying constructible sheaf complex.
The starting point is the following result \cite{Sa0}[sec.2.4.4]:

\begin{thm}[M. Saito]
 Let $M$ be a pure $m$-dimensional complex algebraic manifold. Then one has for $\M\in D^bMHM(M)$ the {\em duality} result (for $j\in \Z$)
\begin{equation}\label{duality-classes1}
 Gr_j^F(DR(\DC \M)) \simeq \DC\left(Gr_{-j}^FDR(\M)\right) \in D^b_{coh}(M) \:.
\end{equation}
Here $\DC$ on the left side is the duality of mixed Hodge modules, wheres $\DC$ on the right hand side is the
{\em Grothendieck duality}
$$\DC=Rhom(\cdot, \omega_M[m]):\: D^b_{coh}(M) \to D^b_{coh}(M)\:,$$
with $\omega_M=\Omega^m_M$ the canonical sheaf of $M$. 
\end{thm}

A priori this is a duality for the corresponding analytic (cohomology) sheaves.
Since $\M$ and $DR(\M)$ can be extended to smooth complex algebraic compactification $\bar{M}$, one can apply Serre's GAGA
to get the same result also for the underlying algebraic (cohomology) sheaves. 

\begin{cor}[Characteristic classes and duality]
 Let $Z$ be a complex algebraic variety with {\em dualizing complex} $\omega_Z^{\bullet}\in  D^b_{coh}(Z)$,
so that the {\em Grothendieck duality transformation} $\DC=Rhom(\cdot ,\omega_Z^{\bullet})$ induces a duality involution
$$\DC: G_0(Z)\to G_0(Z) \:.$$
Extend this to $G_0(Z)\otimes\Z[y^{\pm 1}]$ by $y\mapsto 1/y$. Then the
motivic Hodge Chern class transformation $MHC_y$ commutes with duality $\DC$:
\begin{equation}\label{duality-classes2}
MHC_y(\DC(\cdot))= \DC(MHC_y(\cdot)): K_0(MHM(Z))\to G_0(Z)\otimes\Z[y^{\pm 1}] \:.
\end{equation}
\end{cor}

Note that for $Z=pt$ a point this reduces to the duality formula
(\ref{chi-y-dual}) for the $\chi_y$-genus. For dualizing complexes and (relative) Grothendieck duality we refer to \cite{Ha, Con, LH} as well as
\cite{FM}[Part I,sec. 7]). Note that for $M$ smooth of pure dimension $m$, one has
$$\omega_M[m] \simeq \omega_M^{\bullet} \in  D^b_{coh}(M) \:.$$
Moreover, for a proper morphism $f: X\to Z$ of complex algebraic varieties one has the relative Grothendieck duality isomorphism
$$Rf_*\left(Rhom(\FC ,\omega_X^{\bullet})\right)  \simeq  Rhom(Rf_*\FC ,\omega_Z^{\bullet}) \quad \text{for $\FC \in D^b_{coh}(X)$,}$$
so that the duality involution 
$$\DC: G_0(Z)\otimes\Z[y^{\pm 1}]\to G_0(Z)\otimes\Z[y^{\pm 1}]$$
commutes with proper push down. Since $K_0(MHM(Z))$ is generated by classes $f_*[\M]$, with $f: M\to Z$ proper morphism from a pure dimensional complex algebraic manifold $M$ (and $\M\in MHM(M)$),
it is enough to prove (\ref{duality-classes2}) in the case $Z=M$ a  pure dimensional complex algebraic manifold, in which case it directly follows
from Saito's result (\ref{duality-classes1}). \\

For a systematic study of the behaviour of the Grothendieck duality transformation $\DC:  G_0(Z)\to G_0(Z)$ with respect to exterior products and smooth pullback,
we refer e.g. to \cite{FL} and  \cite{FM}[Part I,sec. 7], where a corresponding ``bivariant'' result is stated. Here we  only point out that the  dualities $(\cdot)^{\vee}$ and $\DC$
commute with the {\em pairings}
of corollary \ref{pairing-class}:
\begin{equation}\begin{split}
 f^*\left((\cdot)^{\vee}\right) &\cap (\DC(\cdot))= \DC\left( f^*(\cdot) \cap (\cdot) \right)
:\\ K^0_{alg}(N)\otimes \Z[y^{\pm 1}] &\times G_0(Z)\otimes \Z[y^{\pm 1}]  \to G_0(Z)\otimes \Z[y^{\pm 1}] \:,
\end{split}\end{equation}
and similarly
\begin{equation}\begin{split}
 f^*\left((\cdot)^{\vee}\right) &\cap (\DC(\cdot))= \DC\left( f^*(\cdot) \cap (\cdot) \right)
:\\ K_0(VmHs^g(N)) &\times K_0(MHM(Z)) \to K_0(MHM(Z))\:.
\end{split}\end{equation}

Here the last equality needs only be checked for classes $[IC_S(\LL)]$, with $S\subset Z$ irreducible of dimension $d$
and $\LL$ a good variation of pure 
Hodge structures on a Zariski dense open smooth subset $U$ of $S$, and $\VB$ a good variation of pure 
Hodge structures on $N$.
But then the claim follows from
$$f^*(\VB)\otimes IC_S(\LL)\simeq IC_S(f^*(\VB)|U\otimes \LL) $$
and (\ref{duality-IC}) in the form
\begin{align*}
{\mathcal D}\left(IC_S(f^*(\mathbb{V})|U\otimes {\mathcal L})\right) 
&\simeq IC_S\left( (f^*(\mathbb{V})|U\otimes {\mathcal L})^{\vee} \right) (d)\\
& \simeq IC_S\left( f^*(\mathbb{V}^{\vee})|U\otimes {\mathcal L}^{\vee} \right) (d) \:.
\end{align*}

\begin{remark}
Also the Todd class transformation $td_*: G_0(\cdot)\to H_*(\cdot)\otimes \Q$ commutes with duality
(compare with \cite{Fu}[ex.18.3.19] and \cite{FM}[Part I, cor.7.2.3]), if the duality
involution $\DC:  H_*(\cdot)\otimes \Q\to  H_*(\cdot)\otimes \Q$ in homology is defined as
$\DC:= (-1)^i\cdot id$ on $H_i(\cdot)\otimes \Q$. So also the unnormalized Hirzebruch class
transformation $MH\tilde{T}_{y*}$ commutes with duality, if this duality in homology
is extended to $H_*(\cdot)\otimes \Q[y^{\pm 1}]$ by ``$\;y\mapsto 1/y$".
\end{remark}

As a final result of this paper, we have the

\begin{prop}\label{MHT=C}
Let $Z$ be a complex algebraic variety, with $[\M]\in K_0(MHM(Z))$. Then
$$MHT_{y*}([\M])\in H_*(Z)\otimes\Q[y^{\pm 1}] \subset 
H_*(Z)\otimes\Q[y^{\pm 1},(1+y)^{-1}]\:,$$
so that the specialization $MHT_{-1*}([\M])\in H_*(Z)\otimes\Q$ for $y=-1$ is defined.
Then
\begin{equation}
MHT_{-1*}([\M]) = c_*([rat(\M)]) =: c_*(\chi_{stalk}([rat(\M)]))
\in H_*(Z)\otimes\Q
\end{equation}
is the rationalized MacPherson Chern class of the underlying constructible sheaf complex
$rat(\M)$ (or the constructible function $\chi_{stalk}([rat(\M)])$). In particular
\begin{equation}
 MHT_{-1*}(\DC[\M]) = MHT_{-1*}([\DC\M]) = MHT_{-1*}([\M]) \:.
\end{equation}
\end{prop}

Here $\chi_{stalk}$ is the transformation form the diagram (\ref{motfunct}).
Similarly, all the transformations from this diagram (\ref{motfunct}), like $\chi_{stalk}$ and $rat$,
commute with duality $\DC$. This implies already the last claim, since $\DC=id$ for algebraically constructible functions
(compare \cite{Sc}[sec.6.0.6]).
So we only need to prove the first part of the proposition.
Since $MHT_{-1*}$ and $c_*$ both commute with proper push down, we can assume
$[\M]=[j_*\VB]$, with $Z=\bar{M}$ a smooth pure dimensional complex
algebraic manifold, $j: M\to \bar{M}$ a Zariski open inclusion with complement $D$ a
normal crossing divisor with smooth irreducible components, and $\VB$ a good
variation of mixed Hodge structures on $M$. So
$$MH\tilde{T}_{y*}([j_*\VB])= ch^*\left( MHC^y(Rj_*L)\right) \cap MH\tilde{T}_{y*}([j_*\Q_M])
\in H_*(\bar{M})\otimes\Q[y^{\pm 1}]$$
by (\ref{MHCj*-fomula}) and the {\em multiplicativity} of the Todd class transformation $td_*$.
Introduce the {\em twisted Chern character}
\begin{equation}
ch^{(1+y)}: K^0_{alg}(\cdot)\otimes \Q[y^{\pm 1}]  \to H^*(\cdot)\otimes\Q[y^{\pm 1}]\;:
\:\:[\VV]\cdot y^j\mapsto \sum_{i\geq 0}\; ch^i([\VV])\cdot (1+y)^i \cdot y^j\:,
\end{equation}
with $ch^i([\VV])\in H^i(\cdot)\otimes\Q$ the $i$-th componenent of $ch^*$. Then
one easily gets
$$MHT_{y*}([j_*\VB])= ch^{(1+y)}\left( MHC^y(Rj_*L)\right) \cap MHT_{y*}([j_*\Q_M])
\in H_*(\bar{M})\otimes\Q[y^{\pm 1},(1+y)^{-1}] \:.$$

But $[j_*\Q_M]=\chi_{Hdg}(j_*[id_M])$ is by (\ref{j-dual}) in the image of 
$$\chi_{Hdg}: M_0(Var/\bar{M})= K_0(Var/\bar{M})[\Lef^{-1}] \to K_0(MHM(\bar{M})) \:.$$
So for $MHT_{y*}([j_*\Q_M])$ we can apply the following special case of proposition \ref{MHT=C}:

\begin{lemma}\label{T=c}
The transformation  
$$T_{y*}=MHT_{y*}\circ \chi_{Hdg}: M_0(Var/Z) \to H_*(Z)\otimes\Q[y^{\pm 1},(1+y)^{-1}]$$
takes values in $H_*(Z)\otimes\Q[y^{\pm 1}] \subset H_*(Z)\otimes\Q[y^{\pm 1},(1+y)^{-1}]$,
with 
$$T_{-1*}=T_{-1*}\circ \DC = c_*\circ can:  M_0(Var/Z) \to H_*(Z)\otimes\Q \:.$$
\end{lemma}

Assuming this lemma, we get from the following commutative diagram, that the specialization $MHT_{-1*}([j_*\VB])$ for $y=-1$ exists:
$$\begin{CD}
 H^*(\cdot)\otimes\Q[y^{\pm 1}] \times H_*(\cdot)\otimes\Q[y^{\pm 1},(1+y)^{-1}]
 @> \cap >> H_*(\cdot)\otimes\Q[y^{\pm 1},(1+y)^{-1}]\\
 @A incl. AA @AA incl. A\\
 H^*(\cdot)\otimes\Q[y^{\pm 1}] \times H_*(\cdot)\otimes\Q[y^{\pm 1}] 
  @> \cap >> H_*(\cdot)\otimes\Q[y^{\pm 1}]\\
 @V y=-1 VV @VV y=-1 V\\
 H^*(\cdot)\otimes\Q \times H_*(\cdot)\otimes\Q 
  @> \cap >> H_*(\cdot)\otimes\Q \:.
\end{CD}$$

Moreover $ch^{(1+y)}\left( MHC^y(Rj_*L)\right)$ specializes for $y=-1$ just to
$$rk(L)=ch^0([\;\overline{\LL}\;]) \in H^0(\bar{M})\otimes \Q \:,$$ 
with $rk(L)$ the rank of the local system $L$ on $M$. 
So we get
$$MHT_{-1*}([j_*\Q_M])=rk(L)\cdot c_*(j_*1_M)= c_*(rk(L)\cdot j_*1_M) \in H_*(\bar{M})\otimes\Q \:,$$
with $rk(L)\cdot j_*1_M =\chi_{stalk}(rat([j_*\VB]))$. \\

It remains to prove the lemma \ref{T=c}. But all transformations $T_{y*}, \DC, c_*$ and $can$ commute with pushdown for proper maps.
Moreover, by resolution of singularities and additivity, $M_0(Var/Z)$ is generated by classes $[f: N\to Z]\cdot \Lef^k$ ($k\in \Z$), with
$N$ smooth pure $n$-dimensional and $f$ proper. So it is enough to prove that
$T_{y*}([id_N]\cdot  \Lef^k)\in H_*(N)\otimes\Q[y^{\pm 1}]$, with
$$T_{y*}([id_N]\cdot  \Lef^k)= T_{y*}\left(\DC ([id_N]\cdot  \Lef^k)\right) = c_*\left(can([id_N]\cdot  \Lef^k)\right) \:.$$

But by the {\em normalization condition} for our characteristic class transformations one has (compare \cite{BSY}):
$$T_{y*}([id_N])= T_y^*(TN)\cap [N] \in  H_*(N)\otimes\Q[y] \:,$$
with $T_{-1*}([id_N])= c^*(TN)\cap [N]= c_*(1_N)$.
Similarly 
$$T_{y*}([\Lef])=\chi_y([\Q(-1)])=-y \quad \text{ and} \quad  can([\Lef])=1_{pt}\:,$$ 
so that
$$T_{y*}([id_N]\cdot  \Lef^k)\in H_*(N)\otimes\Q[y^{\pm 1}]$$
by the multiplicativity of $MHT_{y*}$ for exterior products (with a point space). Moreover
$$T_{-1*}([id_N]\cdot  \Lef^k)= c_*(1_N)=c_*\left( can([id_N]\cdot  \Lef^k)\right) \:.$$

Finally $\DC([id_N]\cdot  \Lef^k)= [id_N]\cdot \Lef^{k-n}$ by definition of $\DC$, so that
$$T_{-1*}([id_N]\cdot  \Lef^k)= T_{-1*}\left(\DC ([id_N]\cdot  \Lef^k)\right) \:.$$

\section*{Acknowledgements}
This paper is an extended version of an expository  talk  given at the workshop
``Topology of stratified spaces'' at MSRI Berkeley in September 2008. Here I would like to thank the
organizers (G. Friedman, E. Hunsicker, A. Libgober and L. Maxim) for the invitation to this workshop.
 I also would like to thank S. Cappell, L. Maxim and S. Yokura for some discussions on the subject of this paper.

\providecommand{\bysame}{\leavevmode\hbox
to3em{\hrulefill}\thinspace}


\begin{thebibliography}{10}

\bibitem{At} M. F. Atiyah, \emph{The signature of fiber bundles}, in
\emph{Global Analysis (Papers in Honor of K. Kodaira)}, 73--84,
Univ. Tokyo Press, Tokyo, 1969.

\bibitem{Ba} M. Banagl, \emph{Topological Invariants of Stratified
Spaces}, Springer Monographs in Mathematics, Springer Verlag
Berlin-Heidelberg 2007.

\bibitem{Ba2} M. Banagl, \emph{The signature of singular spaces and its refinements to generalized homology theories},
home-page of the author,
in \emph{Topology of Stratified Spaces}, MSRI Publications 58, 2011 (to appear).


\bibitem{BCS} M. Banagl, S.E. Cappell, J.L. Shaneson, \emph{Computing twisted signatures and $L$-classes of stratified spaces},  Math. Ann.  326 (2003), 589-623.


\bibitem{BFM} P. Baum, W. Fulton, R. MacPherson, \emph{Riemann-Roch for singular varieties}, Publ. Math. I.H.E.S. 45, 101-145 (1975).

\bibitem{BFM2} P. Baum, W. Fulton, R. MacPherson, \emph{Riemann-Roch and topological K-theory for singular varieties},
  Acta Math. 143 (1979), 155--192 

\bibitem{BBD}  A.A. Beilinson,  J. Bernstein, P. Deligne,  \emph{Faisceaux pervers}, Ast\'{e}risque 100
(1982).




\bibitem{Bi} F.Bittner,  \emph{The universal Euler characteristic for varieties of characteristic zero},
  Comp. Math. 140 (2004), 1011-1032. 

\bibitem{BSY} J.P. Brasselet, J. Sch\"urmann, S. Yokura, \emph{Hirzebruch classes and motivic Chern classes of singular
spaces}, Journal of Topology and Analysis 2, (2010), 1--55.






\bibitem{CMS0} S.E. Cappell, L.G. Maxim, J.L. Shaneson,
\emph{Euler characteristics of algebraic varieties},
Comm. in Pure and Applied Math. 61 (2008), 409--421.


\bibitem{CMS} S.E. Cappell, L.G. Maxim, J.L. Shaneson,
\emph{Hodge genera of algebraic varieties, I}, 
Comm. in Pure and Applied Math. 61 (2008), 422--449.

\bibitem{CLMS} S.E. Cappell, A. Libgober, L.G. Maxim, J.L. Shaneson,
\emph{Hodge genera of algebraic varieties, II},  
 Math. Annalen 345 (2009), 925--972.

\bibitem{CLMS2} S.E. Cappell, A. Libgober, L.G. Maxim, J.L. Shaneson,
\emph{Hodge genera and characteristic classes of complex algebraic
varieties},  Electron. Res. Announc. Math. Sci. 15 (2008), 1--7.


\bibitem{CS} S.E. Cappell, J.L. Shaneson,  \emph{Stratifiable maps and topological invariants},  J. Amer. Math. Soc.  4  (1991),  no. 3, 521--551

\bibitem{CKS} E. Cattani, A. Kaplan, W. Schmid, \emph{$L^2$ and intersection cohomologies for a polarizable variation of Hodge structure},
Inv. Math. 87 (1987), 217--252.

\bibitem{CHS} S.S. Chern, F. Hirzebruch, J.-P. Serre, \emph{On the index of a fibered
manifold},  Proc. Amer. Math. Soc.  8  (1957), 587--596.


\bibitem{Con} B. Conrad, \emph{Grothendieck Duality and Base Change},
Lecture Notes in Mathematics, Vol. 1750. Springer-Verlag, 2000.

\bibitem{De2} P. Deligne, \emph{Th\'{e}or\`{e}me de Lefschetz et crit\`{e}res de d\'{e}g\'{e}n\'{e}rescence de suites spectrales},
Publ. Math. IHES 35 (1968), 107--126.


\bibitem{De} P. Deligne, \emph{Equation diff\'erentielles a point
singular r\'egulier}, Springer 1969.

\bibitem{De1} P. Deligne, \emph{ Th\'{e}orie des Hodge II},
   Publ. Math. IHES 40 (1971), 5--58.

\bibitem{De3} P. Deligne, \emph{Th\'{e}orie des Hodge III},
   Publ. Math. IHES 44 (1974), 5--78. 

\bibitem{Fu} W. Fulton, \emph{Intersection Theory}, Springer-Verlag, 1981.

\bibitem{FL} W. Fulton, L. Lang, \emph{Riemann-Roch Algebra},
  Springer (1985).

\bibitem{FM} W. Fulton, R. MacPherson, \emph{Categorical framework for the study of  
  singular spaces},
  Memoirs of the AMS 243 (1981).

\bibitem{GM2} M. Goresky, R. MacPherson, \emph{Intersection Homology II}, Invent. Math., 71 (1983), 77-129.

\bibitem{Ha} R. Hartshorne,  \emph{Residues and duality},
Lecture Notes in Mathematics, Vol. 20. Springer-Verlag, 1966.

\bibitem{H} F. Hirzebruch,  \emph{Topological methods in algebraic geometry}, Springer, 1966.


\bibitem{Ka} M. Kashiwara, \emph{A study of a variation of mixed Hodge structures},  Publ. RIMS 22 (1986), 991--1024.

\bibitem{KK} M. Kashiwara, T. Kawai, \emph{The Poincar\'{e} lemma for variations of polarized Hodge structures},  Publ. RIMS 23 (1987), 345--407.

\bibitem{KS} M. Kashiwara, P. Schapira, \emph{Sheaves on Manifolds},
Springer-Verlag, Berlin, Heidelberg, 1990.

\bibitem{Ke} G. Kennedy, \emph{MacPherson's Chern classes of singular varieties},
Comm. in Algebra 18 (1990), 2821--2839.


\bibitem{Lib} A. Libgober, \emph{Elliptic genera, real algebraic varieties and quasi-Jacobi forms},
arXiv:0904.1026, in \emph{Topology of Stratified Spaces}, MSRI Publications 58, 2011 (to appear).

\bibitem{LH} J. Lipmann, M. Hashimoto, \emph{Foundations of Grothendieck Duality for Diagrams of Schemes},
Lecture Notes in Mathematics, Vol. 1960. Springer-Verlag, 2009.

\bibitem{M} R. MacPherson, \emph{Chern classes for singular algebraic varieties},
Ann. of Math. (2)  100  (1974), 423--432.

\bibitem{MSc} L. Maxim, J. Sch\"{u}rmann, \emph{Hodge-theoretic Atiyah-Meyer formulae and the stratified 
   multiplicative property},
   in ``Singularities I: Algebraic and Analytic Aspects'',
   Contemp. Math. 474  (2008), 145-167.


\bibitem{Mey} W. Meyer, \emph{Die Signatur von lokalen Koeffizientensystemen
und Faserb\"undeln}, Bonner Mathematische Schriften 53,
(Universit\"at Bonn), 1972.

\bibitem{P} C. Peters, \emph{Motivic aspects of Hodge theory},
Lecture Notes of the Tata Institute of Fundamental Research, 2010.

\bibitem{PS} C. Peters, J. Steenbrink, \emph{Mixed Hodge structures},
Ergenisse der Mathematik und ihrer Grenzgebiete, 3.Folge Vol. 52, Spriner, 2008.

\bibitem{Sa0} M. Saito, \emph{Modules de Hodge polarisables}, Publ. RIMS 24 (1988), 849--995.

\bibitem{Sa1} M. Saito, \emph{Mixed Hodge Modules}, Publ. RIMS 26 (1990), 221--333.


\bibitem{Sa2} M. Saito, \emph{Introduction to mixed Hodge modules}, Actes du Colloque de Th\'{e}orie de Hodge (Luminy, 1987), Ast\'{e}risque No. 179--180 (1989), 10, 145--162.

\bibitem{Sa3} M. Saito, \emph{On Kollar's Conjecture}, Proceedings of Symposia in Pure Mathematics 52, Part 2 (1991), 509--517.

\bibitem{Sa4} M. Saito, \emph{On the formalism of mixed sheaves},
arXiv:math/0611597.

\bibitem{Sa5} M. Saito, \emph{Mixed Hodge complexes on algebraic varieties}, Math. Ann. 316, 283--331 (2000).

\bibitem{Sab} C. Sabbah, \emph{Hodge theory, singularities and D-modules}, prepint (2007), home-page of the author.

\bibitem{Sch} W. Schmid, \emph{Variation of Hodge structures: the singularities of the period mapping},
Inv. Math. 22 (1973), 211--319.

\bibitem{Sc} J. Sch\"urmann, \emph{Topology of singular spaces and constructible sheaves}, Monografie
Matematyczne, 63. Birkh\"auser Verlag, Basel, 2003.

\bibitem{SY} J. Sch\"urmann, S. Yokura,  \emph{A survey of characteristic classes of singular spaces},
in ``Singularity Theory" (ed. by D. Ch\'eniot et al), Dedicated to
Jean Paul Brasselet on his $60$th birthday, Proceedings of the 2005
Marseille Singularity School and Conference, World Scientific, 2007,
865-952.



\bibitem{Si} P.H. Siegel, \emph{Witt spaces: A geometric cycle theory for $KO$-homology at odd primes},
Amer. J. Math. 105 (1983), 1067�-1105.

\bibitem{SZ} J. Steenbrink, S. Zucker, \emph{Variations of mixed Hodge
structures}, Inv. Math. 80 (1983).

\bibitem{To} B. Totaro, \emph{Chern numbers for singular varieties and elliptic homology},
Ann. of Math. (2)  151  (2000),  no. 2, 757--791.

\bibitem{Voi} C. Voisin, \emph{Hodge Theory and Complex Algebaric Geometry, I},
Cambridge studies in advanced mathematics 76. Cambridge University Press, 2002.


\bibitem{Wae} R. Waelder, \emph{Rigidity of differential operators and Chern numbers of singular spaces},
arXiv:0902.4518, in \emph{Topology of Stratified Spaces}, MSRI Publications 58, 2011 (to appear).

\bibitem{Woo} J. Woolf, \emph{Witt groups of sheaves on toplogical spaces},
Comment. Math. Helv. 83 (2008), 289--326.

\bibitem{Y} S. Yokura, \emph{Motivic characteristic classes},
arXiv:0812.4584, in \emph{Topology of Stratified Spaces}, MSRI Publications 58, 2011 (to appear).

\bibitem{Zu} S. Zucker, \emph{Hodge theory with degenerating coefficients: $L_2$-cohomology in the Poincar\'{e} matric},
Ann. Math. 109 (1979), 415--476.

\end{thebibliography}
\end{document}